\RequirePackage{fix-cm}
\documentclass[smallcondensed]{svjour3}     

\smartqed  
\usepackage{amsmath,amssymb}
\usepackage{graphicx}
\usepackage{mathptmx}
\usepackage{algorithm,algorithmic,natbib}
\usepackage{placeins}
\usepackage[labelformat=simple]{subfig}

\usepackage[usenames, dvipsnames]{color}
\usepackage[english]{babel}

\journalname{}

\begin{document}

\title{Structure preserving integration and model order reduction of skew-gradient reaction-diffusion systems}

\titlerunning{Skew-gradient reaction-diffusion systems}        

\author{ B\"{u}lent Karas\"{o}zen \and  Tu\u{g}ba K\"u\c{c}\"ukseyhan \and Murat Uzunca}

\institute{B\"{u}lent Karas\"{o}zen \at Department of Mathematics \& Institute of Applied Mathematics, Middle East Technical University, Ankara-Turkey\\
              Tel.: +90-312-2105602\\
              Fax: +90-312-2101282\\
              \email{bulent@metu.edu.tr}
              \and
              Tu\u{g}ba K\"u\c{c}\"ukseyhan \at Department of Mathematics \& Institute of Applied Mathematics, Middle East Technical University, Ankara-Turkey
              \and
Murat Uzunca \at Institute of Applied Mathematics, Middle East Technical University, Ankara-Turkey
}

\date{Received: date / Accepted: date}

\maketitle

\begin{abstract}
Activator-inhibitor FitzHugh-Nagumo (FHN) equation is an example for reaction-diffusion equations with skew-gradient structure. We discretize the FHN equation using symmetric interior penalty discontinuous Galerkin (SIPG) method in space and average vector field (AVF) method in time. The AVF method is a geometric integrator, i.e. it preserves the energy of the Hamiltonian systems and energy dissipation of the gradient systems. In this work, we show that the fully discrete energy of the FHN equation satisfies the mini-maximizer property of the continuous energy for the skew-gradient systems. We present numerical results with traveling fronts and pulses for one dimensional, two coupled FHN equations and three coupled FHN equations with one activator and two inhibitors in skew-gradient form.
Turing patterns are computed for fully discretized two dimensional FHN equation in the form of spots and labyrinths.  Because the computation of the Turing patterns is time consuming for different parameters, we applied model order reduction with the proper orthogonal decomposition (POD). The nonlinear term in the reduced equations is computed using the discrete  empirical interpolation (DEIM) with SIPG discretization. Due to the local nature of the discontinuous Galerkin (DG) method, the nonlinear terms can be computed more efficiently than for the continuous finite elements.  The reduced solutions are very close to the fully discretized ones. The efficiency and accuracy of the POD and POD-DEIM reduced solutions are shown for the labyrinth-like patterns.

\keywords{FitzHugh-Nagumo equations \and Gradient systems \and Traveling fronts and pulses \and Turing patterns \and Energy preservation \and Discontinuous Galerkin \and Model order reduction \and Discrete empirical interpolation}
\subclass{35K57 \and  65M60  \and 35B36}
\end{abstract}

\section{Introduction}

The FitzHugh-Nagumo (FHN) equation was proposed for modeling the electrical impulses in a nerve axon as a simplified Hodgkin-Huxley model. Although the FHN equation  is developed as a model in physiology, it is used as a generic model that exhibits many phenomena in excitable or oscillatory chemical media. FHN equation is a singularly perturbed, parameter dependent equation with a complex dynamics. In the literature the most known type of the FHN equation is the one which consists of a partial differential equation (PDE) for the activator and an ordinary differential equation (ODE) for the inhibitor. In this work we consider the so called diffusive FHN equation consisting of two or three PDEs with an activator, and one or two inhibitors. Depending on the relationship between the reaction and diffusion parameters, and the bistable nonlinearity, the specific patterns in one and two dimensional FHN equation are like traveling fronts, pulses, spots or labyrinth like patterns. Energy of the system plays an important role in the selection of these patterns. As a reaction  diffusion equation FHN equation has been used to study the wave dynamics in excitable media and bifurcation phenomena.

For large class of reaction-diffusion systems there exist Lyapunov functions, so that all the solutions converge to equilibria. But it is difficult to derive gradient structure for reaction-diffusion systems, because the reaction terms contain quite general nonlinearities. For certain classes, gradient structures can be constructed. Among these are the activator-inhibitor type reaction-diffusion systems which are not order preserving and their linearized forms around the steady states are not self-adjoint. Yanigada \citep{Yanigada02sps} introduced the concept of skew-gradient systems to investigate the stability of these systems. It was shown in \citep{Yanigada02sps} under some restrictions on the parameters the diffusive FitzHugh Nagumo system and Gierer-Meinhardt system exhibit skew-gradient structure. In short, these activator-inhibitor equations consist of  two gradient systems.

In the following we consider the coupled reaction-diffusion system of the form
\begin{equation}\label{rdesg}
\begin{aligned}
\tau_1 u_t &= d_1 \Delta u + f(u,v)\\
\tau_2 v_t &= d_2 \Delta v + g(u,v)
\end{aligned}
\end{equation}
in a  smooth bounded domain $\Omega$ in $\mathbb{R}^{n} (n \leq 2)$ and on a time period $[0,T]$ with $u=u(x,t)$ and $v=v(x,t)$, $x\in\Omega$, $t\in [0,T]$. Here $\tau_1,\; \tau_2$ and $d_1,\; d_2$ correspond to time scales and diffusion coefficients of $u$ and $v$, respectively.
The system (\ref{rdesg}) has a skew-gradient structure if  there holds \citep{Yanigada02sps}
\begin{equation}\label{skew}
\frac{\partial}{\partial v} \left(\frac{f(u, v)}{\tau_1}  \right) = - \frac{\partial}{\partial u} \left(\frac{g(u, v)}{\tau_2}  \right).
\end{equation}

A particular example of the skew-gradient system is the diffusive FHN equation modeling the transmission of electrical impulses in a nerve axon. For this equation, we can write $f(u, v)$ and $g(u, v)$ in equation \eqref{rdesg} as a summation of separate functions of $u$ and $v$ as
\begin{equation}\label{fhn}
\begin{array}{rll}\tau_1 u_t &=& d_1 \Delta u + f_1(u)+f_2(v) \\
\tau_2 v_t &=& d_2 \Delta v + g_1(u)+g_2(v)
\end{array}
\end{equation}
where $f_1(u)=u(u-\beta)(1-u)$, $f_2(v)=-v + \kappa$, $g_1(u)=u$, $g_2(v)= - \gamma v + \epsilon$ and $\gamma > 0$. The FHN equation \eqref{fhn} is an activator-inhibitor system, where $u$ is the activator since it activates $v$, i.e. leads to an increase of $v$. On the other hand, $v$ is the inhibitor since it leads to a decrease in $u$ and $v$. Moreover, FHN equation (\ref{fhn}) is a skew-gradient system (\ref{rdesg}) through the energy functional
$$E(u, v)=\int_{\Omega} \left(\frac{d_1}{2} |\nabla u|^2 - \frac{d_2}{2} |\nabla v|^2 + F(u, v) \right) \mathrm{d}x$$
with the corresponding potential function
\begin{equation}\label{lyp}
F(u,v)=-\frac{u^4}{4}+\frac{(1+\beta )u^3}{3}-\frac{\beta u^2}{2}-uv+\frac{\gamma v^2}{2}-\epsilon v,
\end{equation}
where $\nabla$ denotes the gradient operator. The first equation of \eqref{fhn} is a gradient flow with the potential $F(u, v)$,  and the second with the potential $-F(u, v)$ \citep{Yanigada02mfrd}. The energy functional $E(u,v)$ does not correspond to the Lyapunov functional, because it is not necessarily non-increasing or non-decreasing in time, but the steady state $(u,v)= (\phi, \psi)$ of a skew-gradient system is stable if it is a mini-maximizer of $E(u,v)$ \citep{Yanigada02sps}.

The stationary uniform solutions of \eqref{fhn} can be found by  solving $f_1(\phi)+f_2(\psi)=0=g_1(\phi)+g_2(\psi)$. For one dimensional skew-gradient systems, the number of solutions of this equation determines the type of the solutions of the FHN equation \eqref{fhn}. If it has only one intersection point, then it is called monostable case. However, if there is three intersection points with two stable and one unstable fixed points, then it is called bistable case. For the monostable case, the solution of FHN equation \eqref{fhn} exhibits traveling pulses and for the bistable case, it exhibits traveling fronts. For two dimensional systems, under some conditions on the parameters, Turing patterns may exist, when the spatially homogeneous steady-state  $(u,v)$ is not a mini-maximizer of $E(u,v)$ \citep{Yanigada02sps}.

In this paper, we use symmetric interior penalty discontinuous Galerkin finite elements (SIPG) \citep{arnold82ipf,riviere08dgm}  for space discretization and  average vector field (AVF) method
 \citep{Cellodoni12ped,Hairer14}
 for time discretization. Discontinuous Galerkin (DG) method  uses discontinuous polynomials for the discretization of PDEs.
 DG methods can support high order local approximations that can vary non-uniformly over the mesh, by which fronts, pulses and layers can be captured better.
 The AVF method is second order in time and preserves for conservative systems, like Hamiltonian, the energy and for gradient systems the energy dissipation. 
 We show that the AVF integrator combined with SIPG space discretization preserves the mini-maximizing property of the skew-gradient systems \citep{Yanigada02mfrd} in the discrete form.

FHN equation has front or pulse solutions in one dimensional domains depending on the parameters. We show in numerical simulations the existence of the multiple pulse and front solutions for the two and three component one dimensional FHN equations in the skew-gradient form.
We also show for two the dimensional FHN equation in skew-gradient form the criteria of the existence of patterns like spots and labyrinths.
 The time evolution of the discrete energy shows that the energy remains constant in long term integration, which proves the stability of the solutions obtained by the numerical solutions.

There has been significant development in the efficient implementation and analysis of the model order reduction (MOR) techniques for parametrized PDEs
\citep{Grepl3}.  The aim of the  model order reduction is the construction of a reduced order model (ROM) of small dimension which essentially captures the dynamics of the full order model (FOM). The most known model order reduction technique is the proper orthogonal decomposition (POD). Even though the POD is a very successful MOR technique  for linear problems, for nonlinear problems the computational complexity of the evaluation of the nonlinear reduced model still depends on the dimension of the FOM. In order to reduce the computational cost several methods are developed so that the nonlinear function evaluations are independent of the dimension of the FOM
 and the computational complexity is proportional to the number of ROMs.  Among the best known of them is the discrete empirical interpolation method (DEIM)
 \citep{chaturantabut10nmr} which is adopted from the empirical interpolation method (EIM) \citep{Barrault04}. The DEIM was originally developed for nonlinear functions which depend componentwise on single variables, arising from the finite difference discretization of nonlinear PDEs.
 For the finite element discretization, the nonlinear functions depend on the mesh and on the polynomial degree of the finite elements. Therefore the efficiency of the POD-DEIM can be degraded. A modified version of POD-DEIM is developed in \citep{Tiso13}  using the unassembled finite elements.
 This reduces the number of nonlinear function calls during the online computation, but the size of the nonlinear snapshots are enlarged, which increases the offline computational cost. In \citep{Heinkenschloss14}, the assembled and unassembled finite element POD-DEIM are compared for parametrized steady state PDEs. 
 In case of  DG discretization, each component of the nonlinear functions depends only on the local mesh  in  contrast to the continuous finite element discretization where the nonlinear function depends on multiple components of the state vector. Therefore the number of POD-DEIM function evaluations for DG discretization is comparable with the finite difference discretization.

The paper is organized as follows. In Section~\ref{sec:dg}, the discretization of FHN equation \eqref{fhn} in space by SIPG  is given. In Section~\ref{sec:avf}, the fully discrete formulation of the FHN equation \eqref{fhn} is presented using the  gradient stable AVF time integrator. The Energy analysis of the fully discrete scheme is given in Section~\ref{sec:energy}. In Section~\ref{sec:front},  numerical simulations are shown for the FHN equation exhibiting front and pulse solutions in one dimensional domains.  We present in Section~\ref{sec:3fhn} the numerical results for the three component FHN equation in skew-gradient form with multiple pulse solutions in one dimensional domains. Turing patterns for two-dimensional FHN equation are shown in Section~\ref{sec:turing}. We show  in Section~\ref{sec:rom} the effectiveness and accuracy of the reduced order models for SIPG discretization with application of the POD-DEIM to the nonlinear term. The paper ends with some conclusions.


\section{Space discretization by discontinuous Galerkin method}
\label{sec:dg}
In last twenty years, the DG methods gained an increasing importance for an efficient and accurate solution of PDEs. Although continuous finite elements (FEM) require continuity of the solution along element interfaces, DG does not require continuity of the solution along the interfaces. In contrast
to the stabilized continuous Galerkin finite element methods, DG methods produce stable discretization without the need for extra stabilization strategies, and damp the unphysical
oscillations.  Due to the local structure of DG discretization, DG methods are parallelizable and adaptive meshing techniques can be implemented efficiently. The DG methods combines the best properties of the finite volume and continuous finite elements methods.

In this section, we will describe the DG discretized  semi-discrete (continuous in time) form of the FHN equation \eqref{fhn} with homogeneous Neumann (zero-flux) boundary condition. The classical (continuous) weak solution of \eqref{fhn} solves for $t\in (0,T]$ the variational formulation
\begin{equation}\label{cweak}
\begin{array}{lll}
(\tau_1 u_t, w_1) + a(d_1;u, w_1) -(f_1(u), w_1) - (f_2(v), w_1) &=& 0, \quad  \forall w_1,\\
(\tau_2 v_t, w_2) + a(d_2;v, w_2) -(g_1(u), w_2) - (g_2(v), w_2) &=& 0, \quad  \forall w_2
\end{array}
\end{equation}
with the initial conditions satisfying
$$
(u(0), w_1) = (u_0, w_1), \qquad
(v(0), w_2) = (v_0, w_2),
$$
where $(\cdot, \cdot):=(\cdot, \cdot)_{\Omega}$ is the $L^2$ inner product over the domain $\Omega$, $w_1$ and $w_2$ are the test functions and $a(d;u, w)=(d\nabla u, \nabla w)_{\Omega}$ stands for the classical bilinear form.



Let $\varepsilon_h$ be the disjoint partition of the domain $\Omega$ with elements (triangles) $\{E_i\}_{i=1}^{N_{el}}\in\varepsilon_h$, where $N_{el}$ is the number of elements in the partition. On $\varepsilon_h$, we set the discrete solution and test function space
$$
D_k=D_k(\varepsilon_h):=\{w \in L^2(\Omega): w_E \in \mathbb{P}_k(E)\; \forall E \in \varepsilon_h\},
$$
where $\mathbb{P}_k(E)$ is the space of polynomials of degree at most $k$ on $E\in \varepsilon_h$.
Multiplying  \eqref{fhn} by $w_1$ and $w_2$ and integrating by using Green's theorem over each mesh element, we obtain the following semi-discrete variational formulation: $\forall t \in (0,T]$, find $u_h(t), \; v_h(t)\in D_k$ satisfying
\begin{equation} \label{dgvar}
\begin{array}{lll}
\left(\tau_1\frac{\partial u_h}{\partial t}, w_1\right) + a_h(d_1;u_h, w_1) - (f_1(u_h), w_1) - (f_2( v_h), w_1) &=& 0, \quad \forall w_1\in D_k,\\
\left(\tau_2\frac{\partial v_h}{\partial t}, w_2\right) + a_h(d_2;v_h, w_2) - (g_1(u_h), w_2) - (g_2(v_h), w_2) &=& 0, \quad \forall w_1\in D_k,
\end{array}
\end{equation}
where the SIPG bilinear form $a_h(d;u, w)$ is given by
\begin{align*}
a_h(d;u, w) &= \sum\limits_{E\in \varepsilon_h}^{} \int_E d \nabla u \cdot \nabla w -
\sum\limits_{e\in\Gamma_h^0} \int_e \{d \nabla u\} \cdot [w] \mathrm{d}s   \\
& \quad - \;   \sum\limits_{e\in\Gamma_h^0} \int_e \{d \nabla w\} \cdot [u] \mathrm{d}s +
\sum\limits_{e\in\Gamma_h^0} \frac{\sigma d}{h_e} \int_e  \{\nabla u\} \cdot [w] \mathrm{d}s,
\end{align*}
where $h_e$ denotes the length of an edge $e$, $\Gamma_h^0$ denotes the set of inter-element faces (edges), $[\;]$ and $\{\;\}$ stand for the jump and average operators respectively. 
Introducing the degrees of freedom $N:=N_{loc}\times {N_{el}}$, where $N_{loc}$ denotes the local dimension on each element depending on the polynomial order $k$, the semi-discrete DG solutions of \eqref{dgvar} are of the form
\begin{equation}\label{dgsoln}
u_h(t)=\sum_{i=1}^{N} u_i(t)\phi_i \; , \qquad v_h(t)= \sum_{i=1}^{N} v_i(t)\phi_i,
\end{equation}
where ${\mathbf u}(t):=(u_1(t),\ldots , u_{N}(t))^T$ and ${\mathbf v}(t):=(v_1(t),\ldots , v_{N}(t))^T$ are the vectors of time dependent unknown coefficients of $u_h$ and $v_h$, respectively, and $\phi:=(\phi_1,\ldots , \phi_{N})^T$ is the vector of basis functions. Plugging \eqref{dgsoln} into the scheme \eqref{dgvar} and choosing $w_1=w_2=\phi_i$, $i=1,\cdots,N$,  we get the system of $2\times N$ dimensional ODEs for the unknown vectors ${\mathbf u}$ and ${\mathbf v}$ as
\begin{equation}\label{dgm}
\begin{aligned}
\tau_1 M {\mathbf u}_t+ S_u {\mathbf u} - F_1({\mathbf u}) - F_2({\mathbf v}) &= 0, \\
\tau_2 M {\mathbf v}_t + S_v {\mathbf v} - G_1({\mathbf u}) - G_2({\mathbf u}) &= 0,
\end{aligned}
\end{equation}
where $M,\; S_u,\; S_v \in \mathbb{R}^{N\times N}$ are the mass matrix and  stiffness matrices, respectively, 
and the remaining are the vectors in $\mathbb{R}^{N}$ of the unknowns ${\mathbf u}$ and ${\mathbf v}$.


\section{Time discretization by the average vector field method}
\label{sec:avf}
Energy stable time discretization methods preserve the dissipative structure of the numerical solution of gradient flow equations and skew-gradient systems like the FHN equation. The small values of the diffusion parameters leads to stiff systems after spatial discretizations. Implicit/implicit-explicit
methods are developed since the explicit methods are not suitable for stiff systems and the fully implicit systems require solution of non-linear equations at each time step. In the semi-implicit schemes, the linear stiff part is treated implicitly and the nonlinear part explicitly, so that at each time step a linear system of equations is solved. However, these methods do not preserve the energy dissipation of the system.
Implicit Euler method and average vector field (AVF) method are energy stable time discretization techniques which are robust with small diffusion parameters. The AVF method is the only second order implicit energy
stable method and it preserves energy decreasing property for the gradient systems and for systems with Lyapunov functionals like the FHN equation. In this work, we apply AVF method to solve the system of ordinary differential equations \eqref{dgm} arising from the semi-discretization of the model problem \eqref{fhn}.

We split the time interval $[0, T]$ into $J$ equally-length subintervals $(t_{k-1},t_k]$ with $0 = t_0 < t_1 < \ldots < t_J = T$ with the uniform step-size $\Delta t=t_k - t_{k-1},\; k= 1, 2, \ldots , J$. 
The AVF method  for an arbitrary ODE $\dot{y}=f(y)$ is given by
\begin{equation}\label{avfff}
\frac{y_{n+1}-y_n}{\Delta t} = \int_0^1 f(\xi y_{n+1}+(1-\xi)y_n)\mathrm{d}\xi,
\end{equation}
where $y_n\approx y(t_n)$ for $n=1,\ldots , J$ and $y_0=y(t_0)$. Applying the AVF formulation \eqref{avfff} to the system of ODEs \eqref{dgm}, in space SIPG discretized fully discrete formulation of the equation \eqref{fhn} reads as 
\begin{equation}\label{lavf}
\begin{aligned}
\left(\tau_1 M + \frac{\Delta t}{2} S_u\right) {\mathbf u}_{n+1} &= \left(\tau_1 M - \frac{\Delta t}{2} S_u \right) {\mathbf u}_n + \Delta t \int\limits_{0}^{1} F_1(\xi {\mathbf u}_{n+1}+(1-\xi){\mathbf u}_n) \mathrm{d}\xi \\
 & \quad + \; \Delta t \int\limits_{0}^{1} F_2(\xi {\mathbf v}_{n+1}+(1-\xi){\mathbf v}_n) \mathrm{d}\xi, \\
\left(\tau_2 M + \frac{\Delta t}{2} S_v\right) {\mathbf v}_{n+1} &= \left(\tau_2 M - \frac{\Delta t}{2} S_v \right) {\mathbf v}_n + \Delta t \int\limits_{0}^{1} G_1(\xi {\mathbf u}_{n+1}+(1-\xi){\mathbf u}_n) \mathrm{d}\xi \\
& \quad + \; \Delta t \int\limits_{0}^{1} G_2(\xi {\mathbf v}_{n+1}+(1-\xi)  {\mathbf v}_n) \mathrm{d}\xi.
\end{aligned}
\end{equation}
The fully discrete  system of nonlinear equations \eqref{lavf} is solved by  Newton's method on  each time-interval $(t_{n-1},t_n]$, $n=0,1,\ldots , J-1$.

\section{Energy analysis of fully discrete scheme}
\label{sec:energy}
 It was proved in \citep{Yanigada02mfrd} that a steady state of a skew-gradient system is stable if and only if it is a mini-maximizer of the energy functional \eqref{energy}:
\begin{itemize}
\item If $u=\xi$ is a local minimizer of $E(u,\eta)$ and $v=\eta$ is a local maximizer of $E(\xi, v)$, then $(u, v)=(\xi, \eta)$ is a mini-maximizer of $E(u,v)$.
\item
If $E(\xi, \bar{v}) \leq E(\xi, \eta) \leq E(\bar{u}, \eta)$ for any neighborhoods $\bar{u},\; \bar{v}$ of $\xi$ and $\eta$, resp., then we say that $(u, v)=(\xi, \eta)$ is a mini-maximizer of $E(u, v)$.
\end{itemize}
We will examine the relation between the stability of a steady state solution $(\xi, \eta)$ of \eqref{fhn} and the mini-maximizing property of the critical point of the energy functional \eqref{energy} for the SIPG-AVF fully discrete system.


Let $u_n$ and $v_n$ denote the solutions of \eqref{fhn} at $t=t_n$. Steady state solutions $(u_n,v_n)=(\xi(x), \eta(x))$ of \eqref{fhn} satisfies
\begin{equation}\label{stdy}
\begin{aligned}
d_1 \Delta \xi + f_1(\xi) + f_2(\eta) &= 0, \\
d_2 \Delta \eta + g_1(\xi) + g_2(\eta) &= 0.
\end{aligned}
\end{equation}
Then the solution of \eqref{stdy} is a critical point of the following energy functional
\begin{equation} \label{energy}
\begin{aligned}
E(u_{n}, v_{n}) &= \frac{d_1}{2}\left\|\nabla u_{n}\right\|^2_{L^2(\varepsilon_h)} - \frac{d_2}{2}\left\|\nabla v_{n}\right\|^2_{L^2(\varepsilon_h)} + (F(u_n, v_n),1)   \\
 & \quad + \sum_{e\in\Gamma_h^0} \left( -(\{d_1\partial_n u_{n}\}, [u_{n}])_e + \frac{\sigma d_1}{2h_e}([u_{n}],[u_{n}])_e \right)  \\
  & \quad - \sum_{e\in\Gamma_h^0} \left( -(\{d_2\partial_n v_{n}\}, [v_{n}])_e + \frac{\sigma d_2}{2h_e}([v_{n}],[v_{n}])_e \right),
\end{aligned}
\end{equation}
where $F(\cdot , \cdot)$ is the potential function \eqref{lyp}. Fully discrete variational formulation reads
\begin{equation}\label{avf}
\begin{aligned}
\frac{\tau_1}{\Delta t} (u_{n+1} - u_{n}, w_1) & =  - \frac{1}{2}  a_h(d_1; u_{n+1}+u_n,w_1) + \int\limits_{0}^{1} (f_1(\xi u_{n+1}+(1-\xi)u_n), w_1) d\xi \\
& \quad +  \int\limits_{0}^{1} (f_2(\xi v_{n+1}+(1-\xi)v^n), w_1) d\xi  \\
\frac{\tau_2}{\Delta t} (v_{n+1} - v_{n}, w_2) & =  - \frac{1}{2}  a_h(d_2; v_{n+1}+v_n,w_2) + \int\limits_{0}^{1} (g_1(\xi u_{n+1}+(1-\xi)u_n), w_2) d\xi \\
&  \quad  + \int\limits_{0}^{1} (g_2(\xi v_{n+1}+(1-\xi)v_n), w_2) d\xi
\end{aligned}
\end{equation}
By taking  $w_1=u_{n+1}-u_n$ and $w_2=v_{n+1}-v_n$ in \eqref{avf}, and using the identity $(a+b, a-b) = (a^2 - b^2,1) $ and the bilinearity of $a_h$, we get
\begin{equation}\label{useidentity}
\begin{aligned}
\frac{\tau_1}{\Delta t} (u_{n+1} - u_{n}, u_{n+1}-u_n)   & =  -\frac{1}{2}  a_h(d_1; u_{n+1},u_{n+1}) + \frac{1}{2}  a_h(d_1; u_{n},u_{n}) \\
&  \quad  + (f_1(\xi u_{n+1}+(1-\xi)u_n), u_{n+1}-u_n) d\xi \\
&   \quad  + (f_2(\xi v_{n+1}+(1-\xi)v_n), u_{n+1}-u_n) d\xi,  \\
\frac{\tau_2}{\Delta t} (v_{n+1} - v_{n}, v_{n+1}-v_n)   & =  -\frac{1}{2}  a_h(d_2; v_{n+1},v_{n+1}) + \frac{1}{2}  a_h(d_2; v_{n},v_{n})  \\
&  \quad   + (g_1(\xi u_{n+1}+(1-\xi)u_n), v_{n+1}-v_n) d\xi  \\
&  \quad  + (g_2(\xi v_{n+1}+(1-\xi)v_n), v_{n+1}-v_n) d\xi.
\end{aligned}
\end{equation}
Using the Taylor expansion for the terms $f_1,\; f_2,\; g_1,\; g_2$ in \eqref{useidentity},
and subtracting $F(u_n, v_n)$ from $F(u_{n+1}, v_{n+1})$, we get
\begin{equation}\label{subtr}
\begin{aligned}
F(u_{n+1}, v_{n+1})-F(u_n, v_n) &\approx   \frac{\partial F}{\partial u} (\xi u_{n+1}+(1-\xi)u_n, v_{n})(u_{n+1}-u_n) \\
& \quad + \frac{\partial F}{\partial v} (\xi u_{n+1}+(1-\xi)u_n, v_{n})(v_{n+1}-v_n).
\end{aligned}
\end{equation}
Using the derivatives $\frac{\partial F}{\partial u} = f_1(u) + f_2(v)$, $\frac{\partial F}{\partial v} = - (g_1(u)+g_2(v))$, which are the skew-gradient conditions \eqref{skew}, and plugging the identity \eqref{subtr} into the system \eqref{useidentity}, we obtain the relation
\begin{equation}\label{endiff}
E(u_{n+1},v_{n+1}) - E(u_n,v_n) \approx -\frac{\tau_1}{\Delta t} \left\|u_{n+1} - u_{n}\right\|_{L^2(\Omega)}^2 + \frac{\tau_2}{\Delta t} \left\|v_{n+1} - v_{n}\right\|_{L^2(\Omega)}^2.
\end{equation}

Now, at the steady state point $(u_n, v_n)=(\xi, \eta)$, $\frac{\tau_2}{\Delta t} \left\|v_{n+1} - v_{n}\right\|_{L^2(\Omega)}^2 = 0$ if $v$ is fixed to $\eta(x)$,  leading according to the relation \eqref{endiff} to $E(\xi_{n+1},\eta) \leq E(\xi_n,\eta)$. In a similar way, if $u$ is fixed to $\xi(x)$, we obtain $E(\xi,\eta_{n+1}) \geq E(\xi,\eta_n)$. Thus, the first equation of \eqref{fhn} describes a gradient flow with  $E(u, \eta)$ and the second equation of \eqref{fhn} describes a gradient flow with  $-E(\xi, v)$, meaning that $(u_n, v_n)=(\xi, \eta)$ is stable as a steady state of \eqref{fhn} since it is a mini-maximizer of $E(u,v)$ \citep{Yanigada02mfrd}.


\section{Traveling fronts and pulses of FitzHugh-Nagumo equation}
\label{sec:front}

Localized structures like fronts and pulses are most-well known one-dimensional waves in reaction diffusion systems \citep{Chen14sas}. Fronts connect two different state of a reaction-diffusion system with a bistable nonlinearity. Pulses exist far away from the homogeneous equilibrium and results from the balance between the dissipation and nonlinearity. Existence of fronts and pulses for the FHN equation \eqref{fhn} are shown in \citep{Chen14sas}.
In this section we derive conditions for the parameters of the FHN equation \eqref{fhn} exhibiting traveling fronts or traveling pulses
which depend on mono/bi-stability of the homogeneous steady state solutions of \eqref{fhn}
\begin{equation}\label{hs}
f(u)-v= 0,\qquad \gamma u - \delta v + \epsilon = 0.
\end{equation}
Here, we consider the bistable cubic nonlinearity $f(u)=u(u-\beta)(1-u),\text{ with } 0<\beta <\frac{1}{2}$. Eliminating   $v=\frac{\gamma u + \epsilon}{\delta}$  from \eqref{hs} and fixing the  parameters $\beta=\frac{2}{25},\; \gamma=1,\; \epsilon=\frac{7}{10}$, we obtain
\begin{equation} \label{intt}
u^3 - \frac{27}{25} u^2 + \left(\frac{2}{25} + \frac{1}{\gamma} \right) u + \frac{7}{10 \gamma} = 0.
\end{equation}
For the bistable case, the equation \eqref{intt} must have three  distinct reel roots, or equivalently its derivative equation
\begin{equation}\label{dintt}
3u^2 - \frac{54}{25} u + \left(\frac{2}{25} + \frac{1}{\gamma}\right)
\end{equation}
must have two distinct roots. In other words, we need the condition
\begin{equation*}
\Delta = \left(\frac{54}{25}\right)^2 - 4\cdot 3\cdot \left( \frac{2}{25}+\frac{1}{\gamma}\right) = \frac{2316}{625} - \frac{12}{\gamma} > 0.
\end{equation*}
Hence, for $\gamma > 3.2383$, the equation \eqref{fhn} attains the bistable case, as a result, the solution of \eqref{fhn} is traveling fronts. Otherwise it will be monostable and the traveling pulses will occur.

\subsection{Traveling fronts}
\label{ex:front1}

\begin{figure}[htb!]
\centering
\subfloat[]{\includegraphics[scale=0.3]{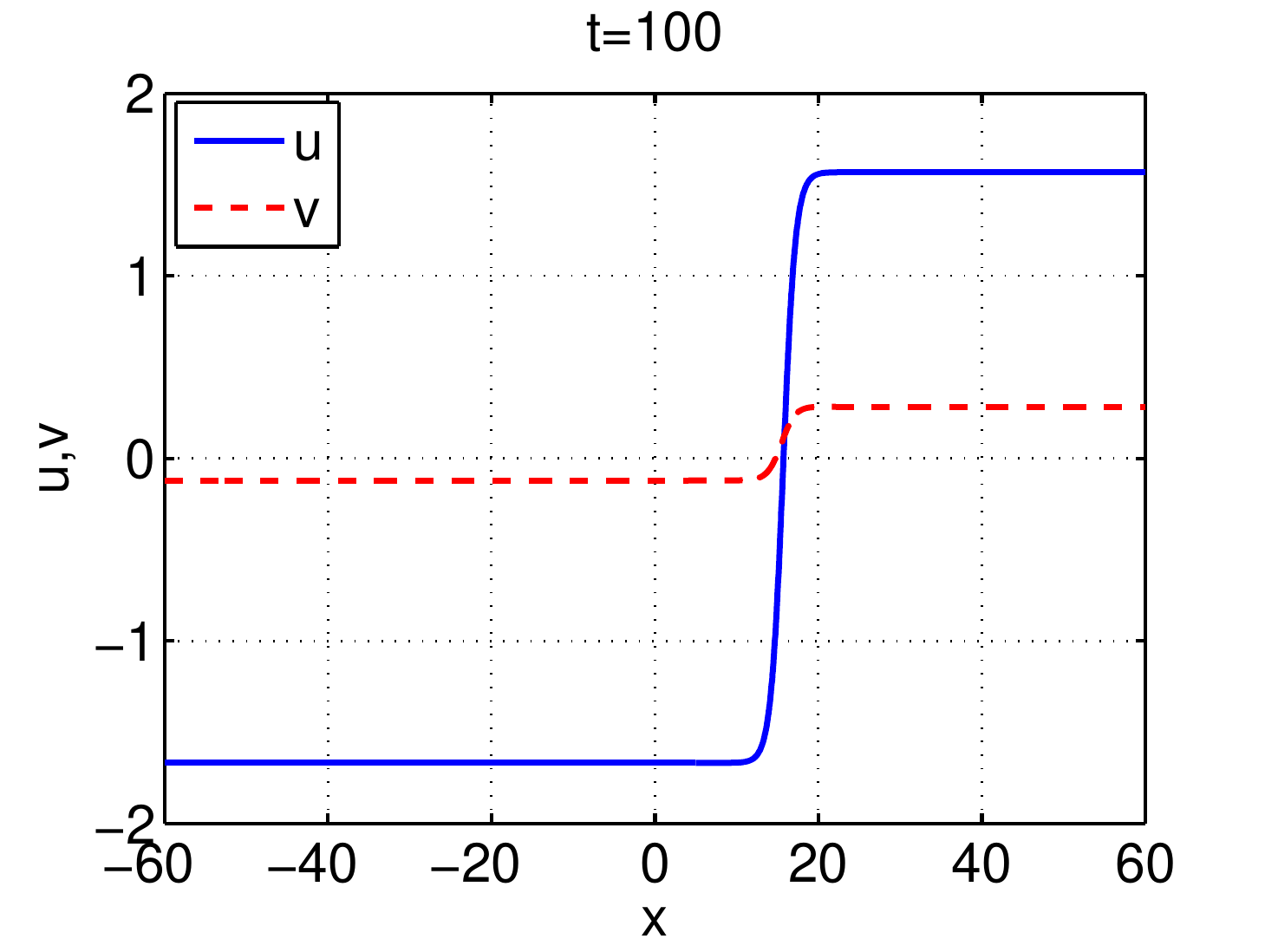}}
\subfloat[]{\includegraphics[scale=0.3]{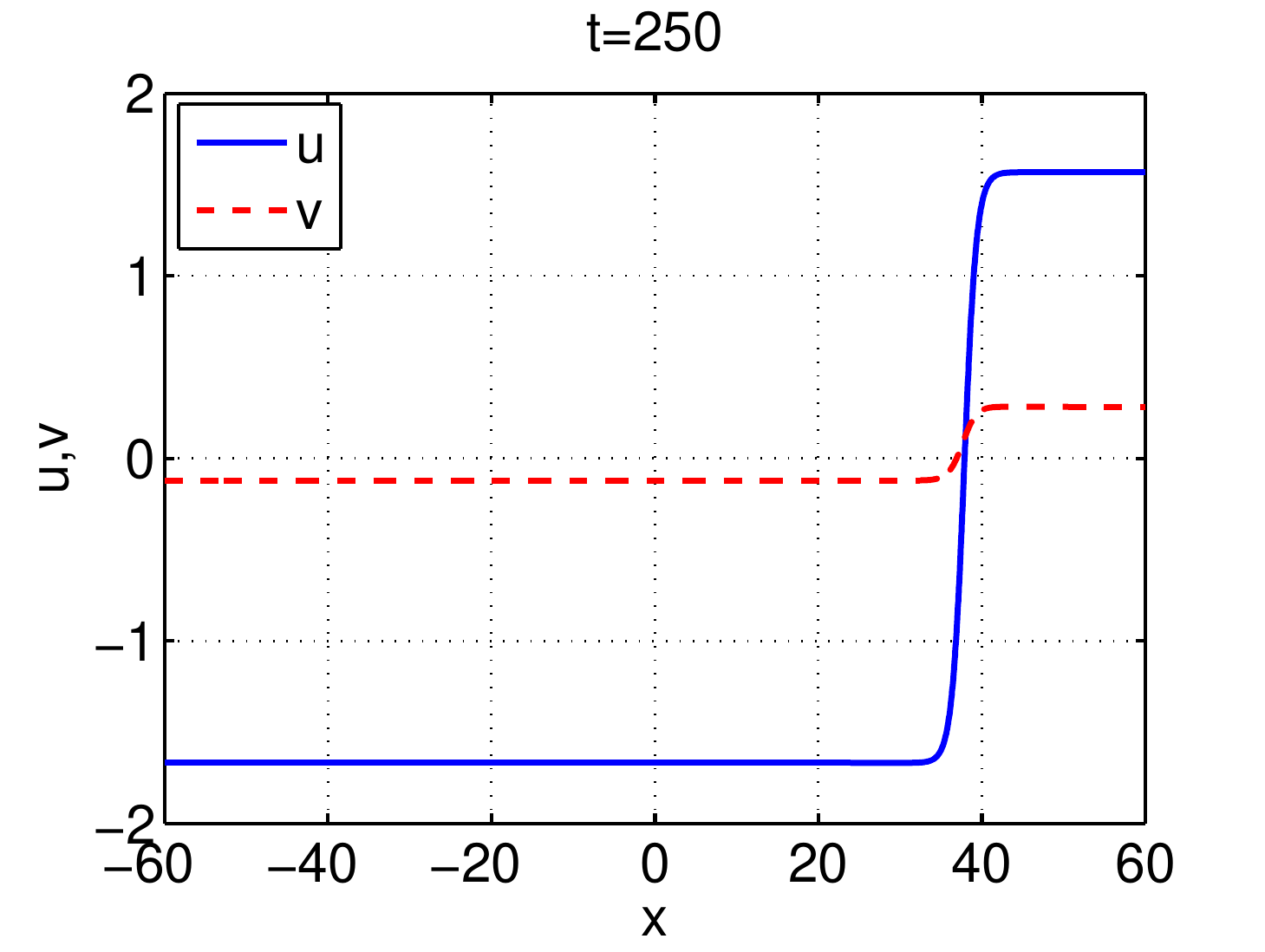}}
\subfloat[]{\includegraphics[scale=0.3]{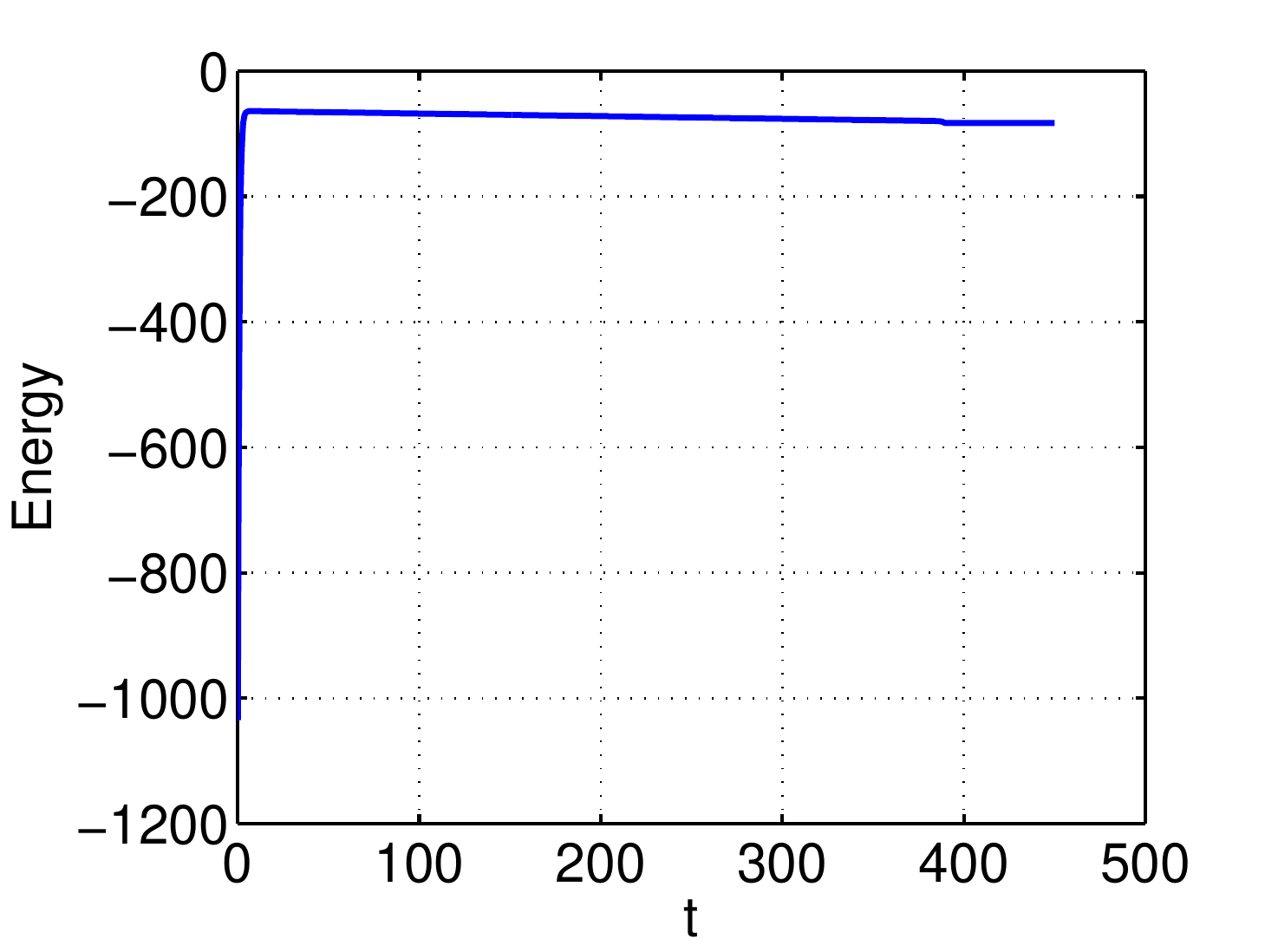}}
\caption{Traveling front solutions and evolution of the discrete energy for Example \ref{ex:front1}}
\label{fig1}
\end{figure}

We consider the FHN equation \eqref{fhn} on $\Omega=[-60,60]$ with the parameters
$d_1=1$, $d_2=1.25$, $\tau=12.5$, $\beta=\frac{1}{3}$, $\gamma=8$ and $\epsilon=0.7$. We set the initial conditions $u(x,0)=\tanh(x)$, $v(x,0)=1-\tanh(x)$, and the spatial mesh size $\Delta x=0.1$ and the temporal step size $\Delta t=0.5$. The traveling front solutions can be clearly seen from Fig.~\ref{fig1}, whereas the discrete energy is decaying very slowly.

\subsection{Traveling pulses}
\label{ex:front2}

\begin{figure}[htb!]
\centering
\subfloat[]{\includegraphics[scale=0.3]{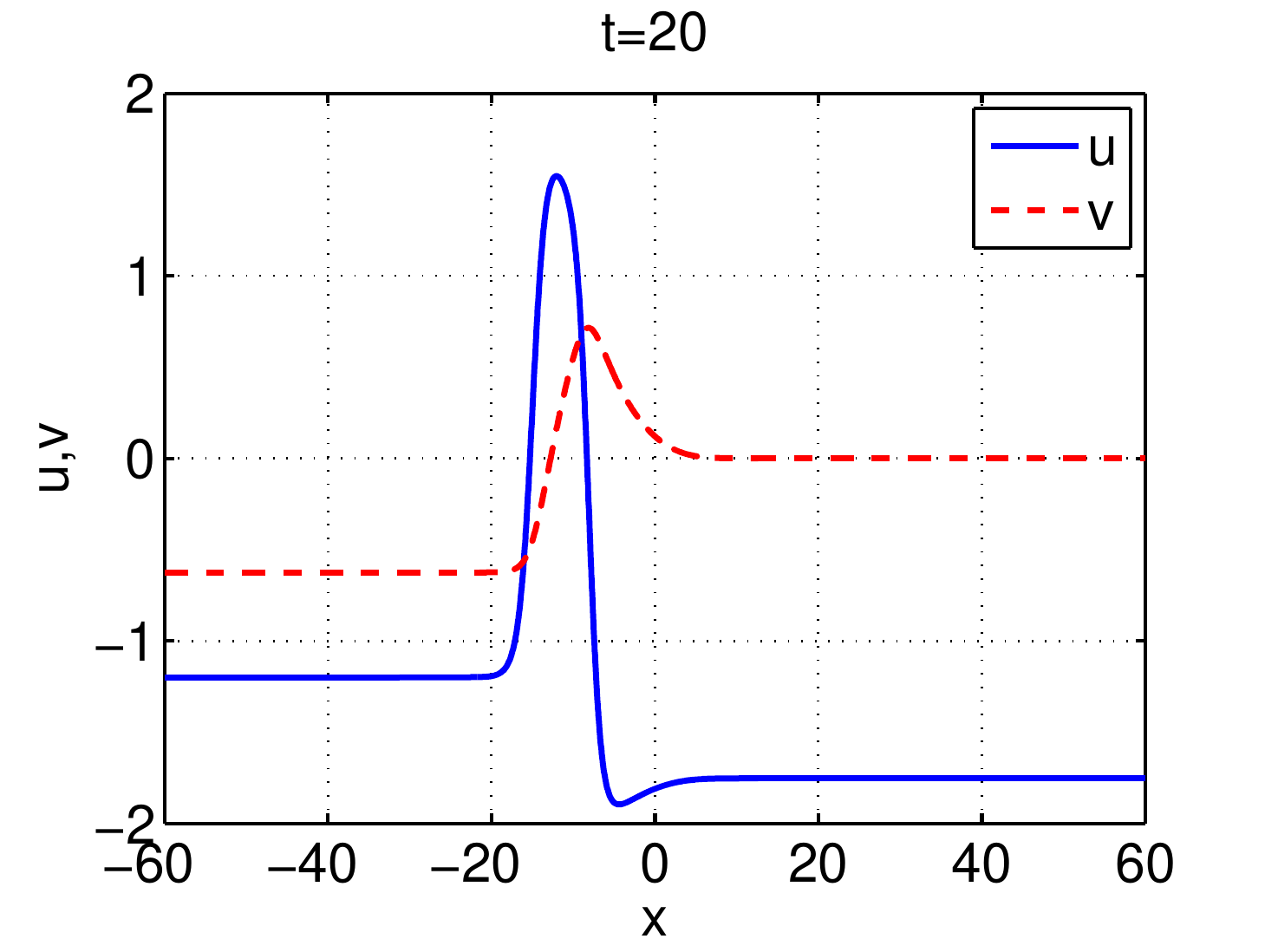}}
\subfloat[]{\includegraphics[scale=0.3]{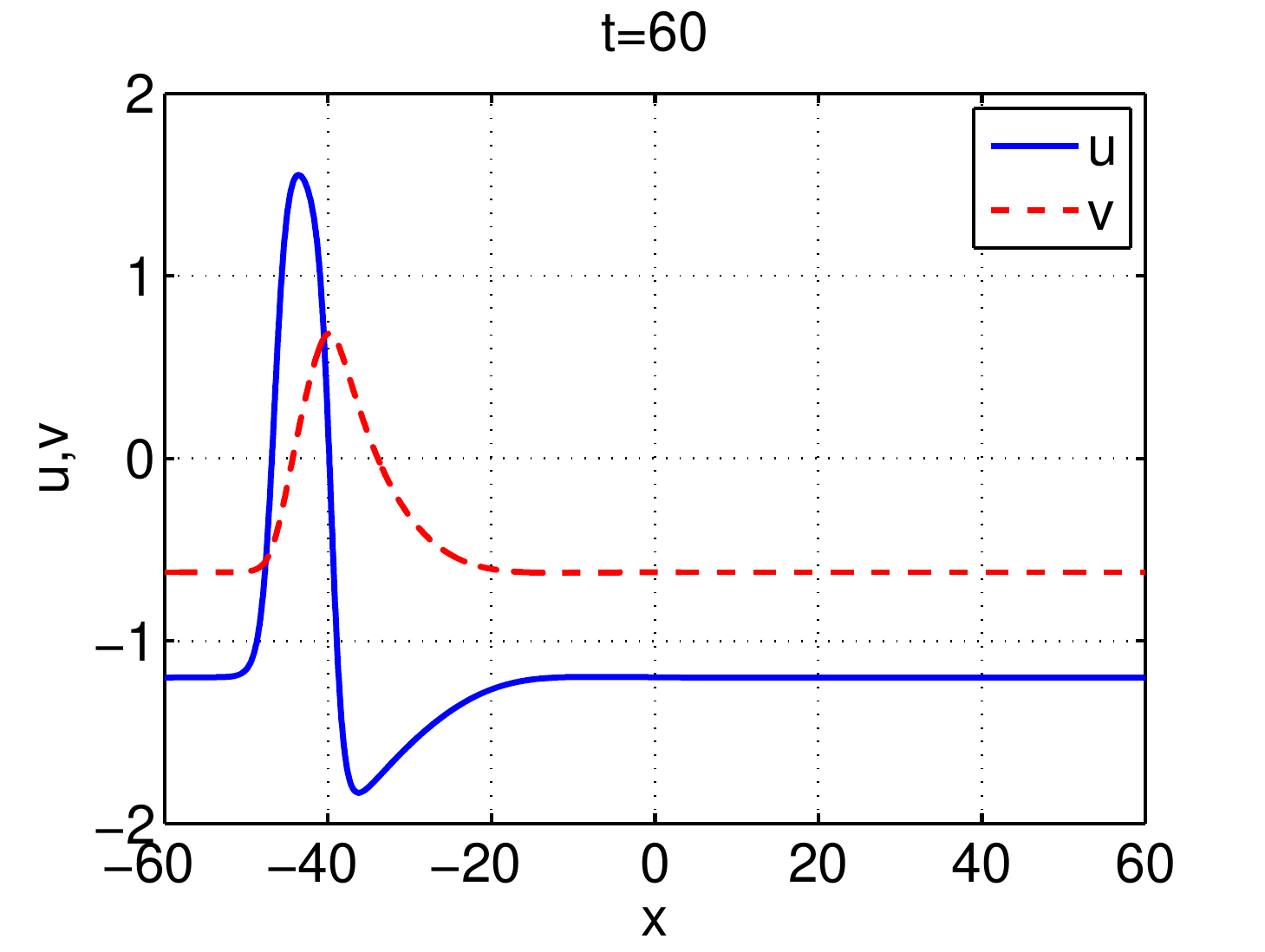}}
\subfloat[]{\includegraphics[scale=0.3]{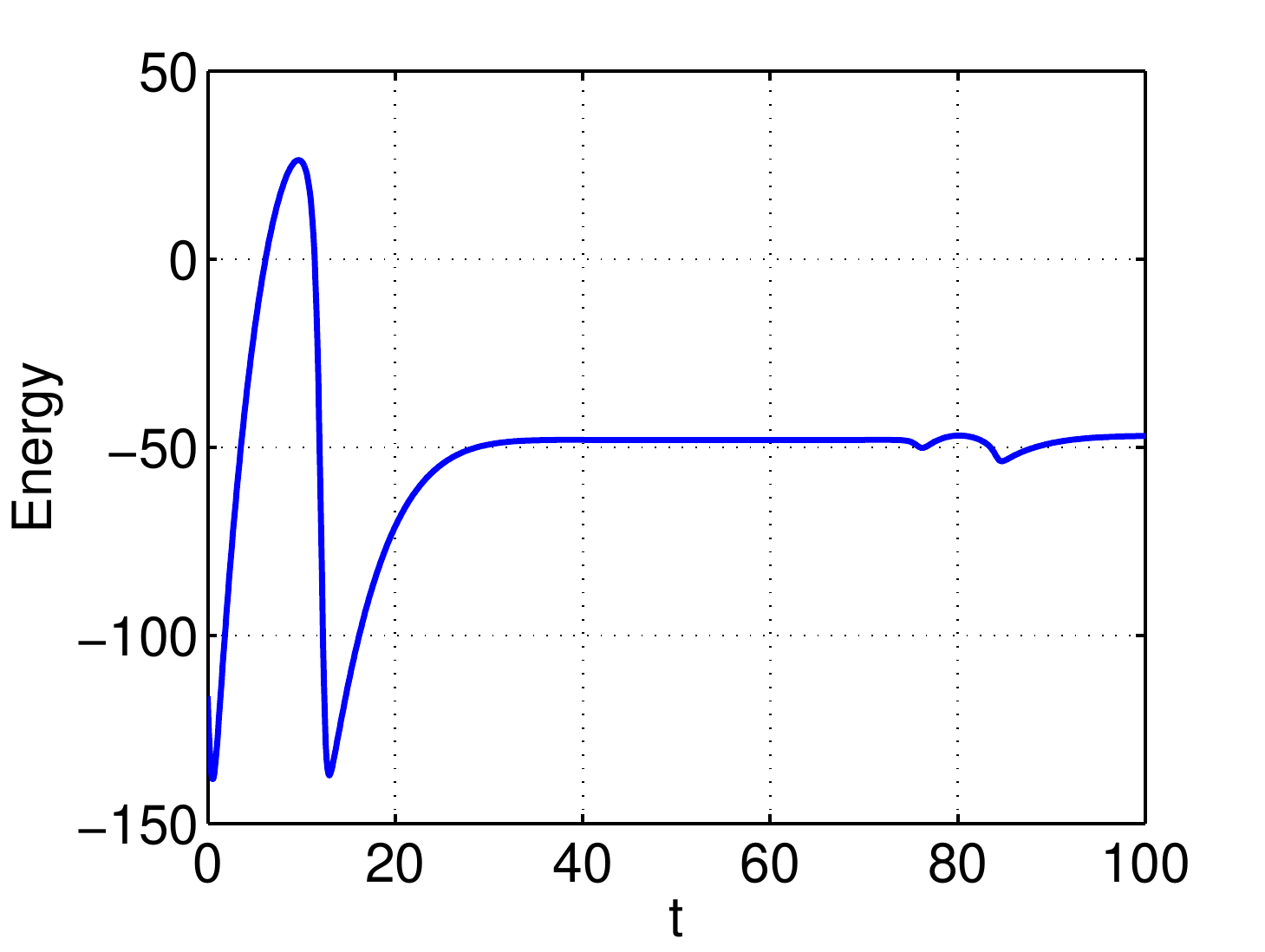}}
\caption{Traveling pulse solutions and evolution of the discrete energy for Example \ref{ex:front2}}
\label{fig2}
\end{figure}

All parameters other than $\gamma=0.8$ are the same as for the traveling front solutions in Example \ref{ex:front1}. Initial conditions are taken as $u(x,0)=\tanh(x)$ and $v(x,0)=-0.6$. The traveling pulse solutions obtained for $\Delta x=0.1$ and $\Delta t=0.1$ are shown in Fig.~\ref{fig2}. The discrete energy remains constant after an oscillation.

\section{Three component FitzHugh-Nagumo system}
\label{sec:3fhn}

The rich dynamics of localized structures like multi-pulse and multi-front solutions were recently investigated for singularly perturbed three component reaction-diffusion systems \citep{guil98,heijster08,heijster11}
\begin{equation} \label{threecomp}
\begin{aligned}
u_t &= u_{xx} + u-u^3-\epsilon(\alpha v+\beta s + \gamma), \\
\tau v_t &= \frac{1}{\epsilon^2} v_{xx} + u - v,  \\
\theta s_t &= \frac{d^2}{\epsilon^2} s_{xx} + u - s,
\end{aligned}
\end{equation}
where $0<\epsilon\ll 1,\; \tau, \theta  > 0,\; d > 1$, and $\alpha$ and $\beta$ denoting the reacting rates and the constant $\gamma$ is the source term. The system \eqref{threecomp} was originally introduced as a model for gas discharge dynamics \citep{guil98}. It became a standard model to study the dynamics and interactions of spatially localized structures like pulses and fronts in one dimensional reaction-diffusion equations. The three component model \eqref{threecomp} is also considered as augmented FHN equation with a second inhibitor $w$ which diffuses more rapidly than the second inhibitor $v$.

The last two equations of \eqref{threecomp} imply  that  all the components of the stationary solutions are equal. Elimination of the components $v$ and $s$ leads to the equation
\begin{equation} \label{hss1}
u^3-u+\epsilon (\alpha u + \beta u + \gamma) =0,
\end{equation}
which has three different roots if its derivative equation $3 u^2 - 1 + \epsilon (\alpha + \beta)$ has two distinct roots
, i.e. $\epsilon (\alpha + \beta) < 1$. In this case, it is called bistable medium since two of the roots of \eqref{hss1} is stable. Hence, \eqref{threecomp} has front solutions, and depending on the initial conditions we can get one-front, two-front or multi-front solutions.

On the other hand, the three component FHN equation \eqref{threecomp} has skew-gradient structure if there hold
\begin{equation} \label{skew3}
\frac{\epsilon \alpha}{\tau} = \frac{1}{\theta}, \quad \frac{\epsilon \beta}{\tau} = \frac{1}{\theta}.
\end{equation}
In addition to the skew-gradient condition \eqref{skew3}, if $\epsilon (\alpha + \beta) > 1$ is satisfied, then we obtain one stationary point. In this situation, the stationary point is always unstable since the Jacobian matrix around this stationary point of \eqref{threecomp} satisfies $\det (J) < 0$ and  oscillatory traveling waves exists  satisfying the skew-gradient condition \eqref{skew3}.

In the following numerical examples we present one pulse (two-front), two pulse (four front) and multi-pulse (multi-front) solutions of \eqref{threecomp} similar to those in \citep{heijster08} for different set of parameters satisfying the skew-gradient condition \eqref{skew3} under the homogeneous Neumann boundary conditions on both ends. In all examples we have taken $\epsilon=0.01$.

\subsection{One pulse solutions}
\label{ex:multi1}

\begin{figure}[htb!]
\centering
\subfloat[]{\includegraphics[scale=0.3]{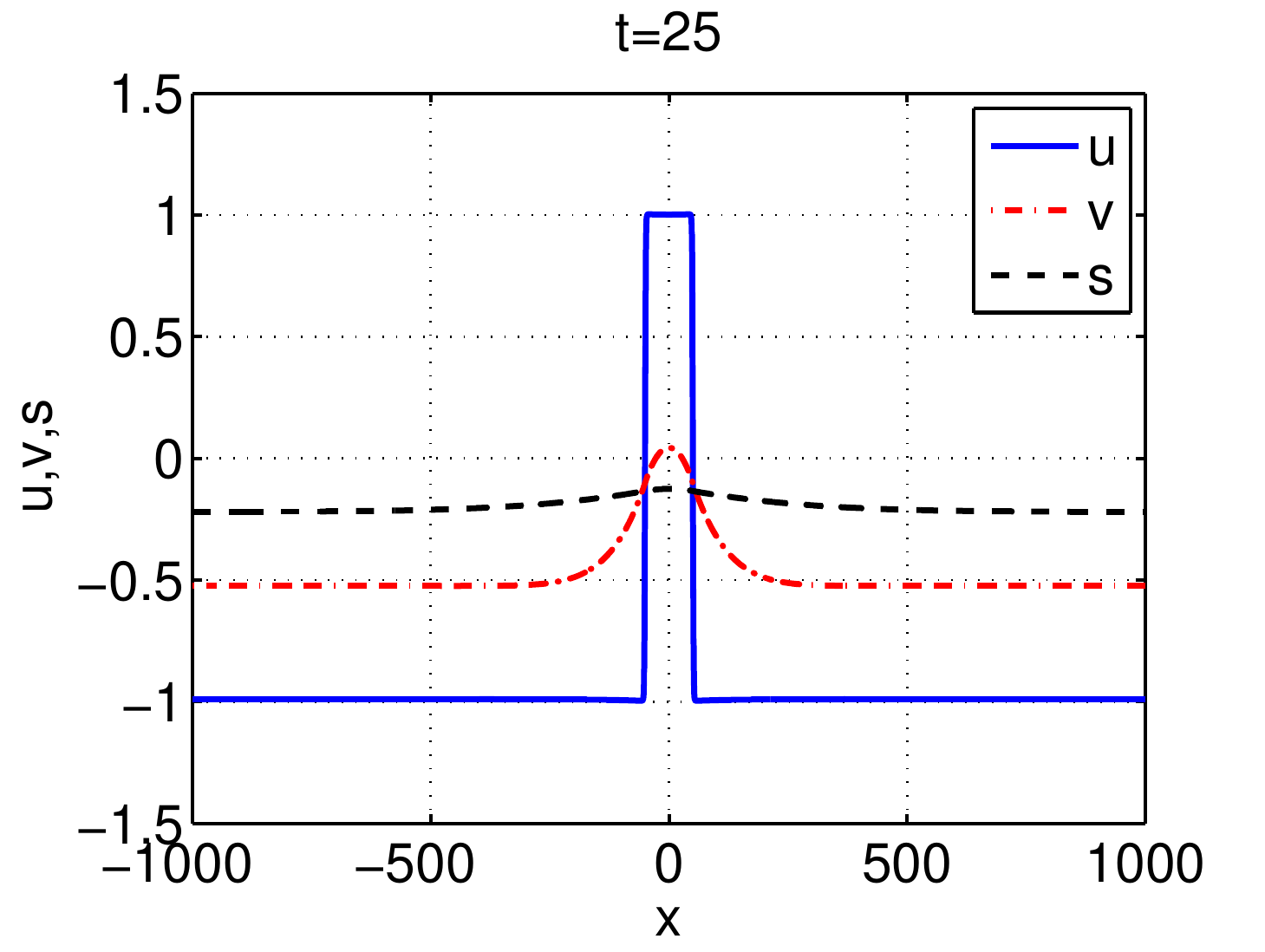}}
\subfloat[]{\includegraphics[scale=0.3]{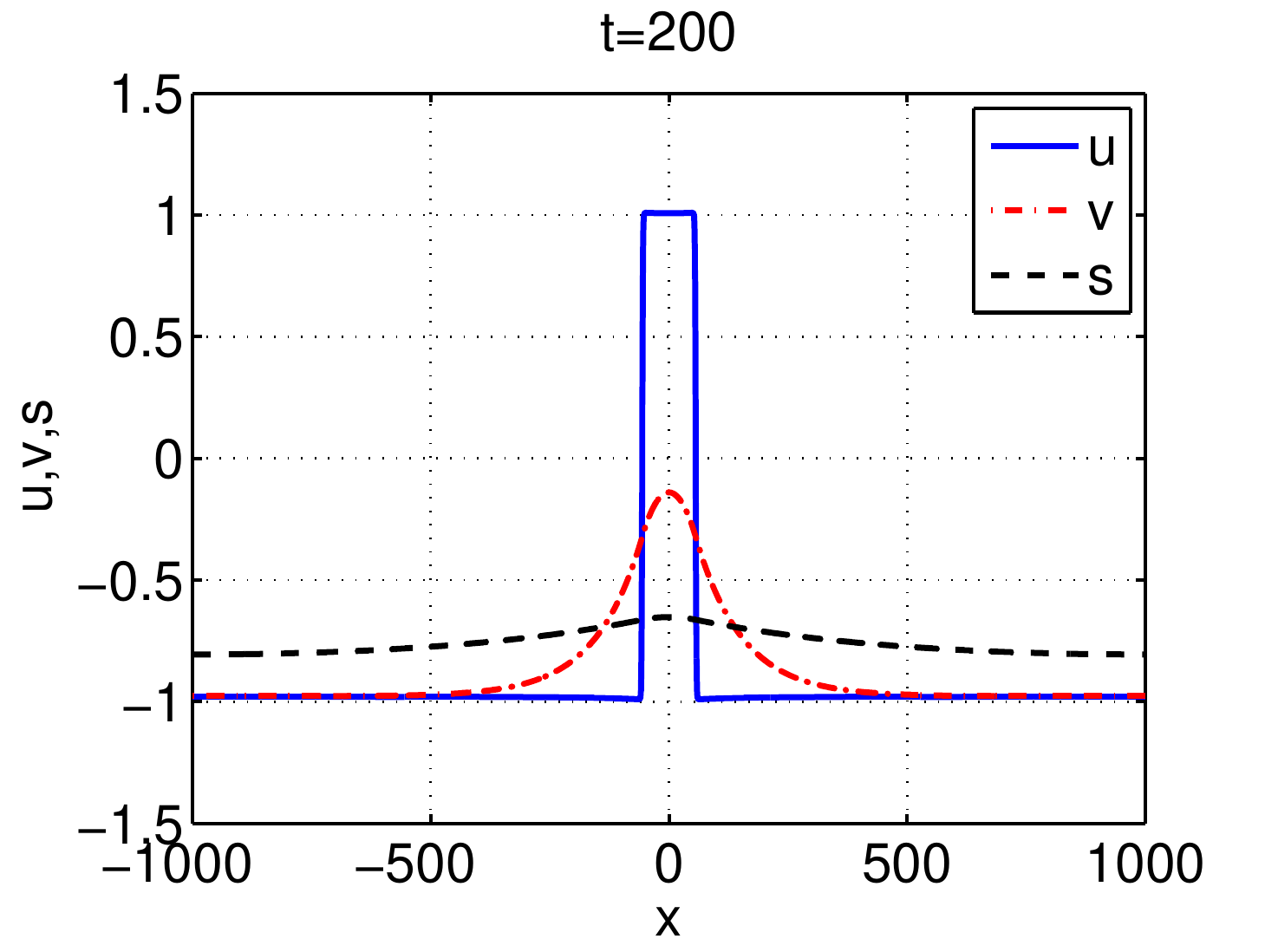}}
\subfloat[]{\includegraphics[scale=0.3]{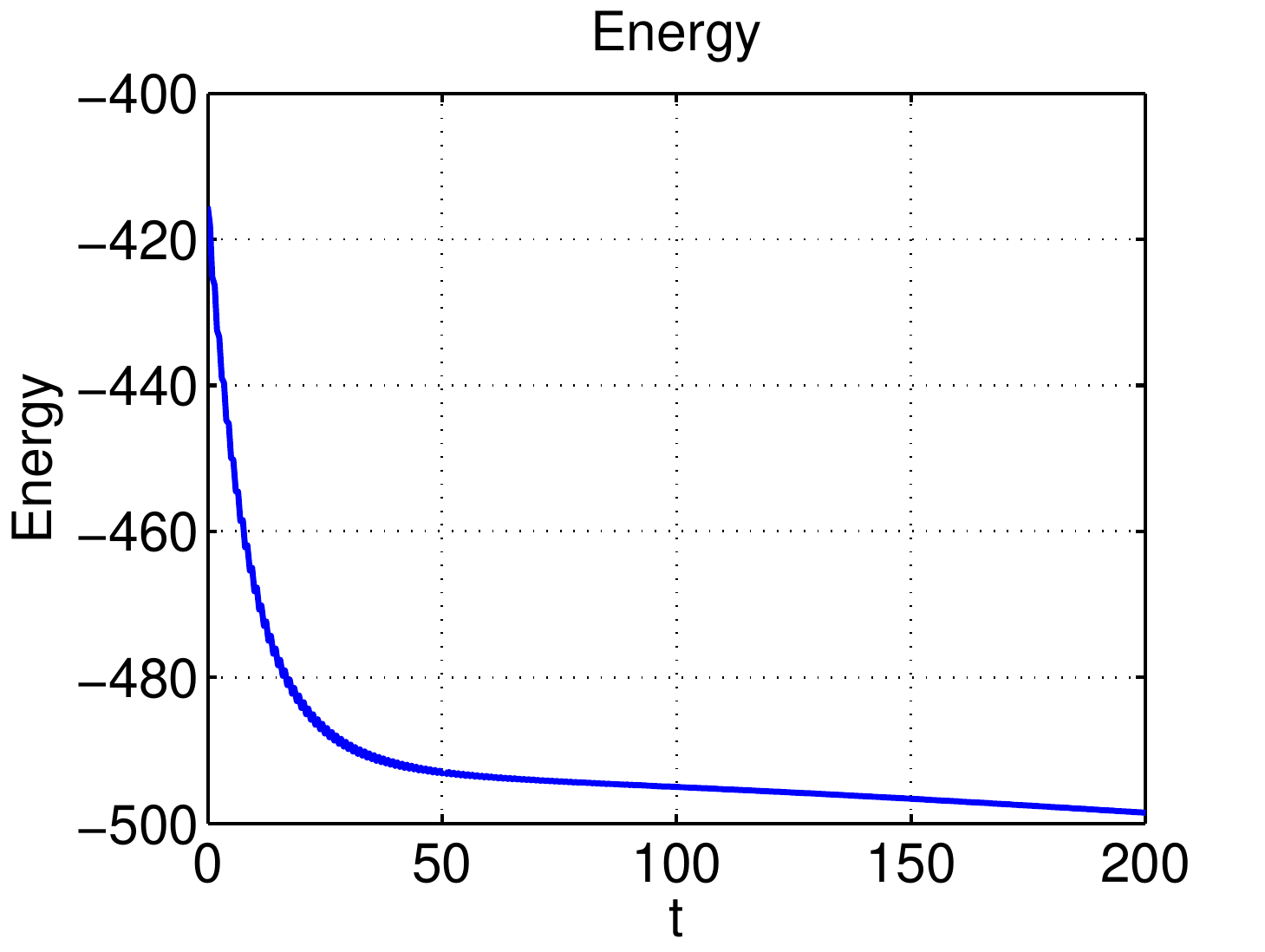}}
\caption{One pulse  solutions and evolution of discrete energy for Example \ref{ex:multi1}}
\label{fig5}
\end{figure}

We consider the space-time domain $x \in [-1000,1000]$ and $t \in [0,200]$ with
$\Delta x= \Delta t = 0.5$, and the initial conditions are given by
\[
u(x,0) = \left\{ \begin{array}{ll}
         1 & \mbox{$x \in [-50,50]$}\\
         -1 & \mbox{otherwise}\end{array} \right. \; , \quad v(x,0)=s(x,0)=0.
\]
We use the set of parameters $(\alpha, \beta, \gamma, d, \tau, \theta)=(3, 1, -0.25, 5, 100/3, 100)$. The solution profiles at $t=25$ and at the final time $t=200$, together with the energy plot, are shown in Fig.~\ref{fig5}.

\subsection{Two pulse solutions}
\label{ex:multi2}

\begin{figure}[htb!]
\centering
\subfloat[]{\includegraphics[scale=0.3]{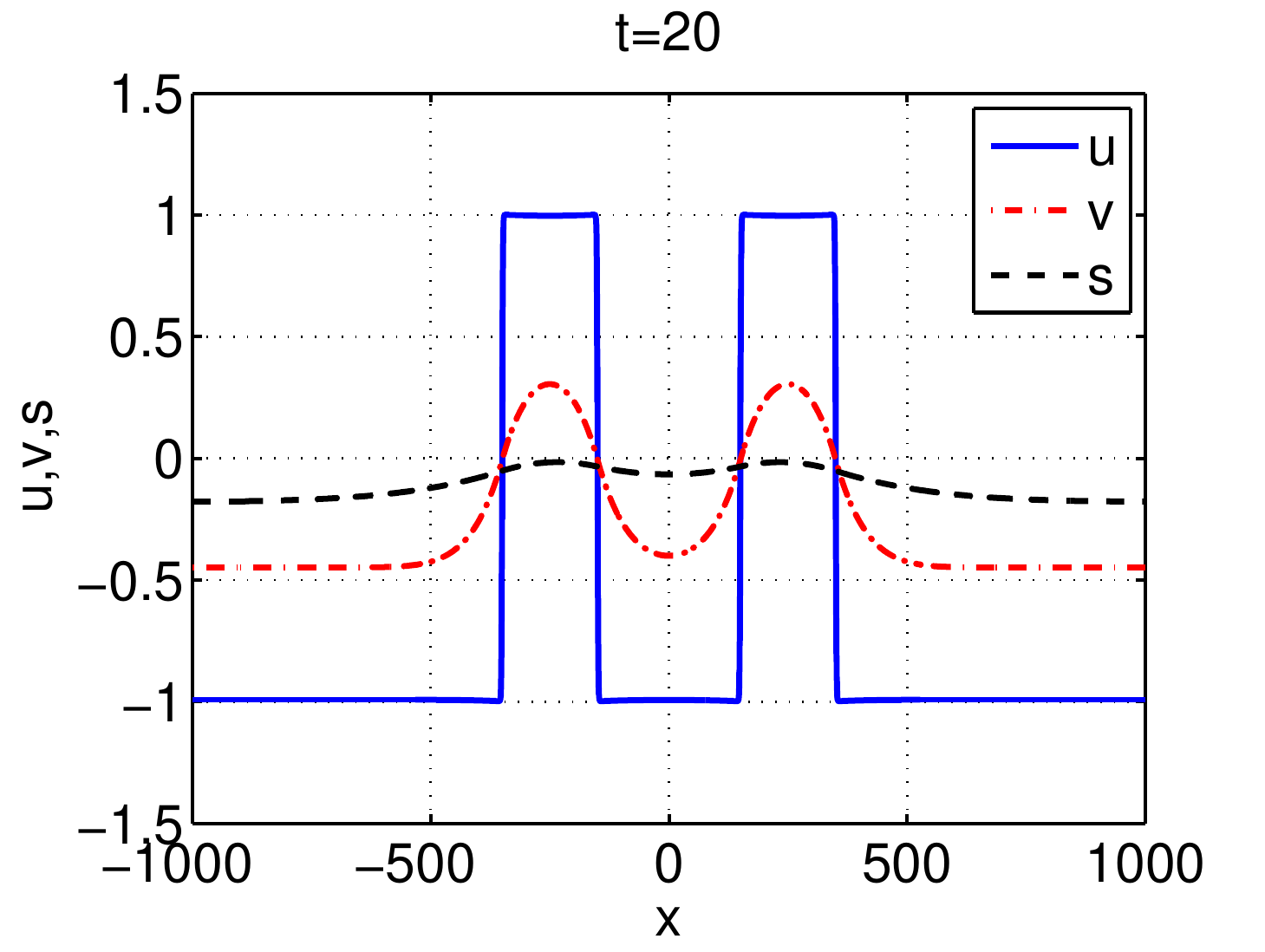}}
\subfloat[]{\includegraphics[scale=0.3]{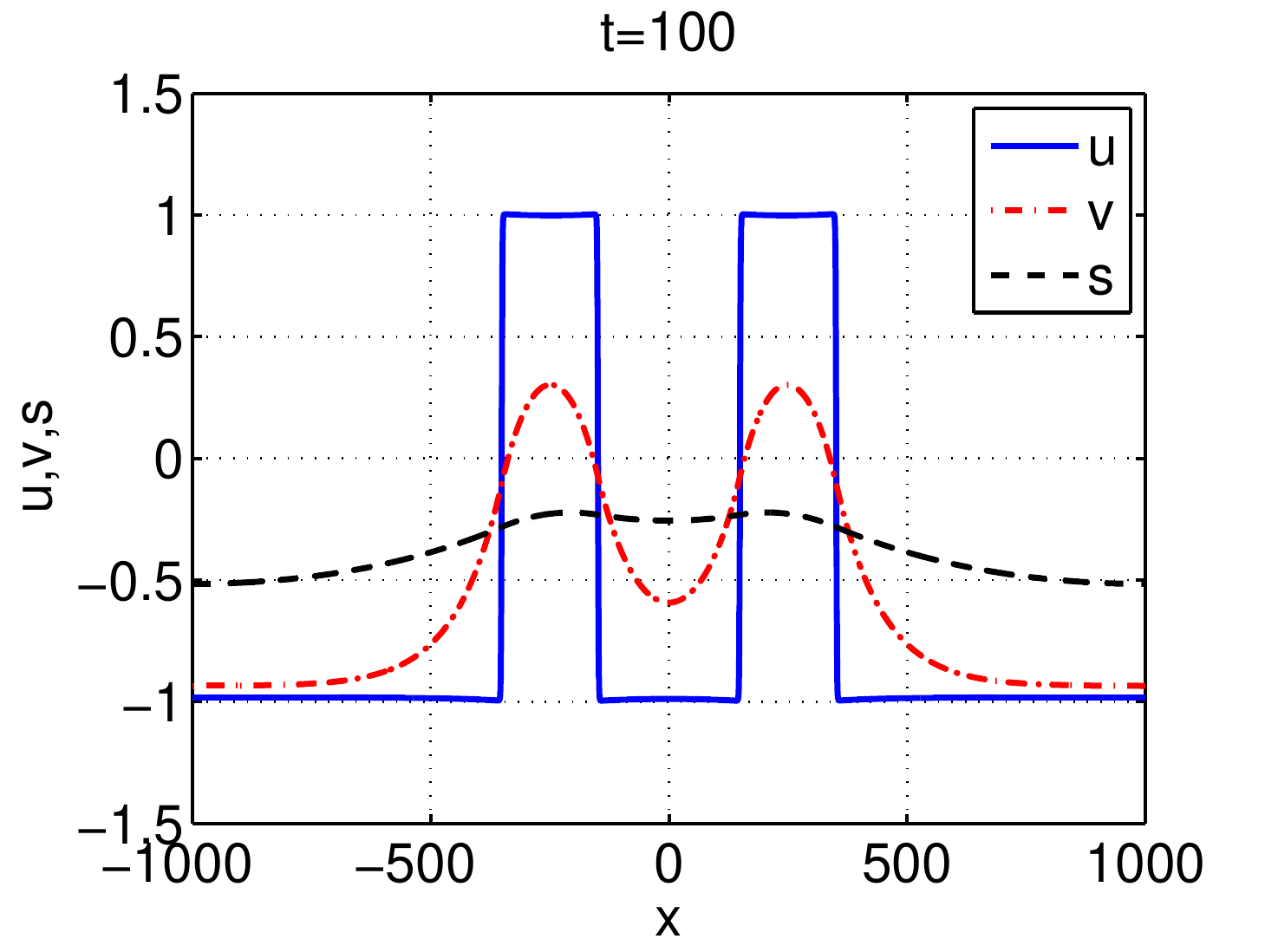}}
\subfloat[]{\includegraphics[scale=0.3]{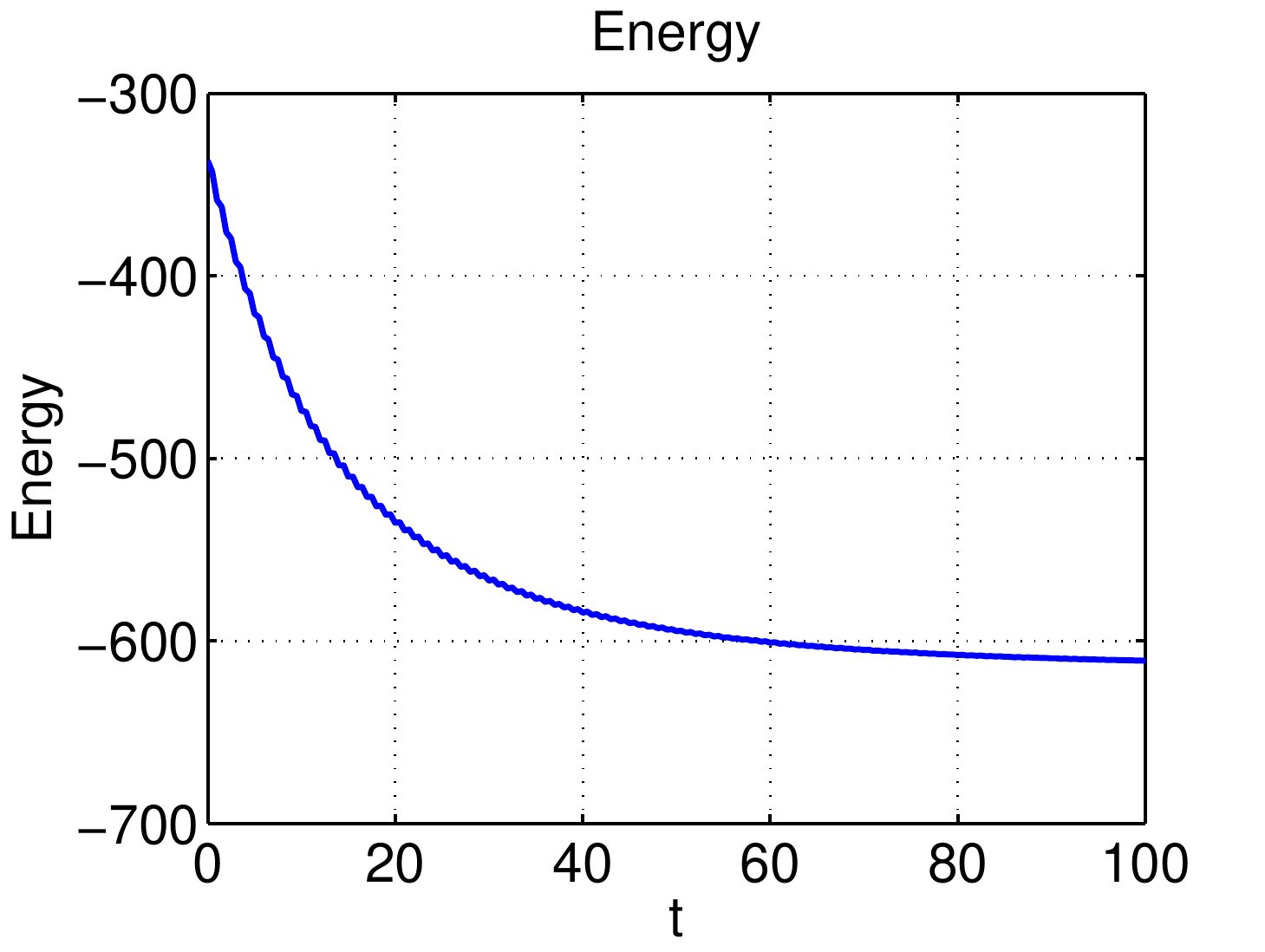}}
\caption{Two pulse solutions and evolution of discrete energy for Example \ref{ex:multi2}\label{fig6}}
\end{figure}

For the two pulse solutions, we use the same settings as in Example \ref{ex:multi1}  with the final time $T=100$. We give the related results in Fig.~\ref{fig6} for the initial conditions
\[
u(x,0) = \left\{ \begin{array}{rcl}
         1 &,& \mbox{if }  x \in [-350, -150]\cup[150, 350] ,\\
          -1 &,& \mbox{otherwise}.\end{array} \right. \; , \quad v(x,0)=s(x,0)=0.
\]
				
\subsection{Multi pulse solutions}
\label{ex:multi3}

\begin{figure}[htb!]
\centering
\subfloat[]{\includegraphics[scale=0.3]{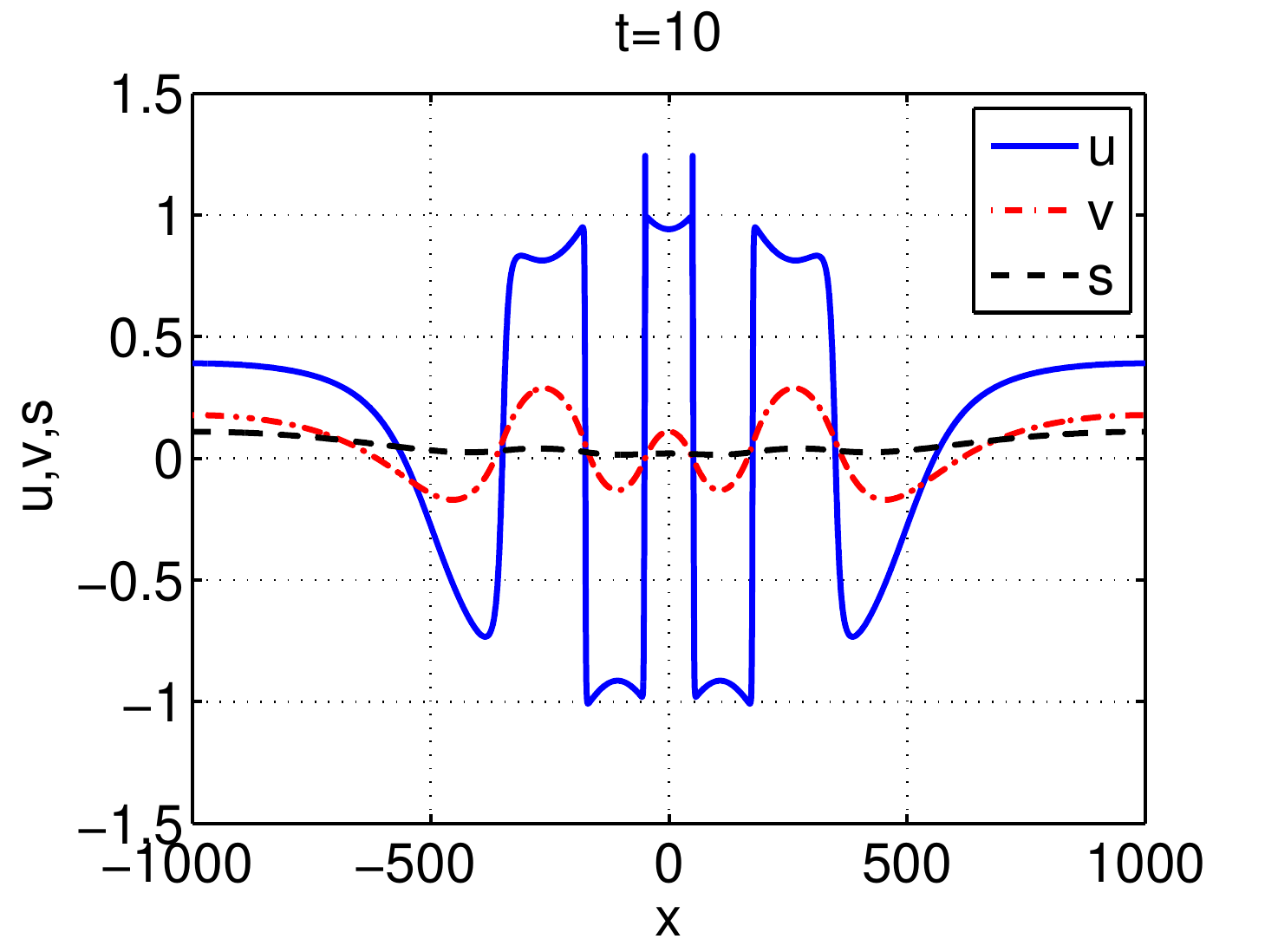}}
\subfloat[]{\includegraphics[scale=0.3]{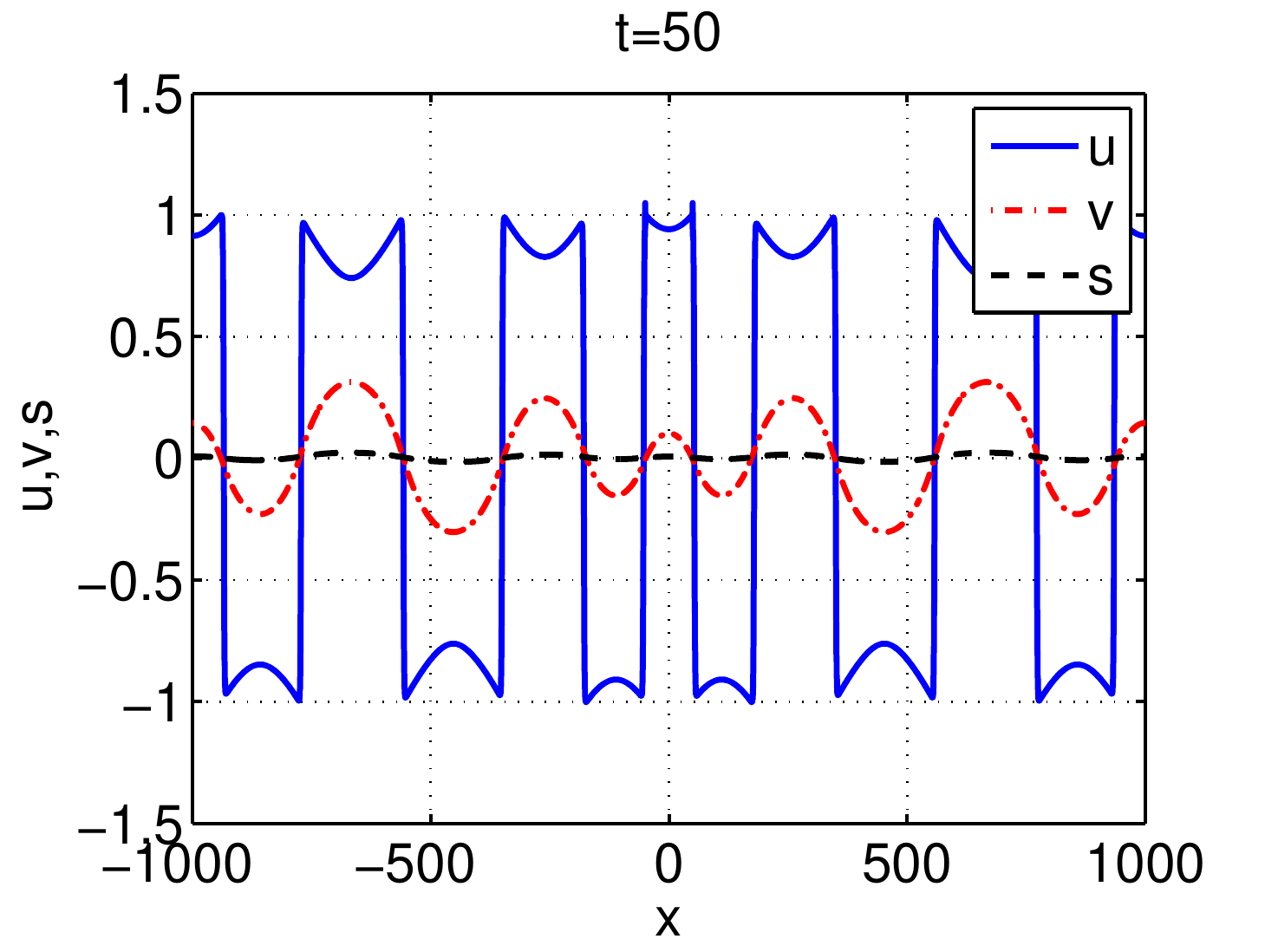}}
\subfloat[]{\includegraphics[scale=0.3]{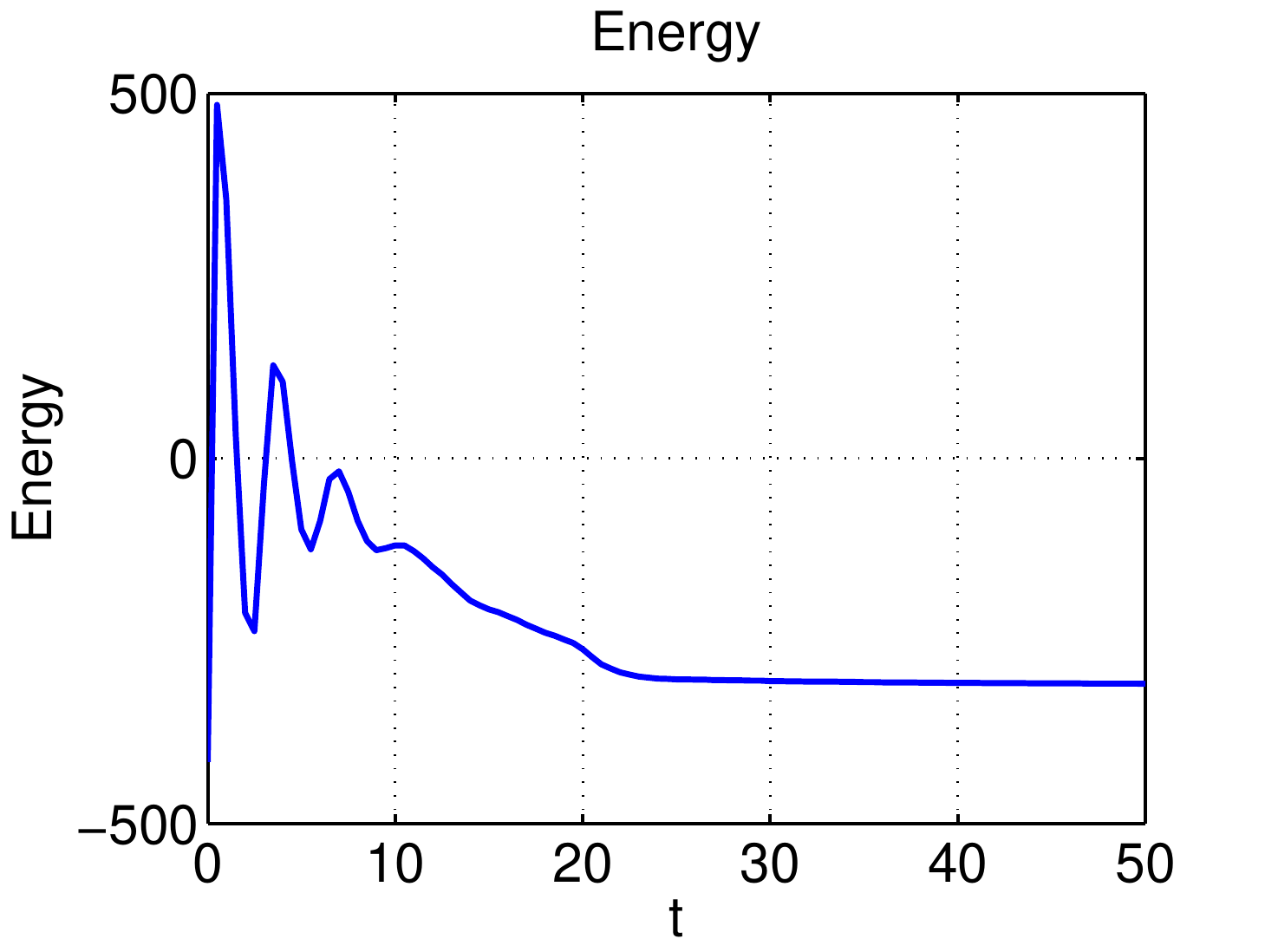}}
\caption{Multi pulse solutions and evolution of discrete energy for Example \ref{ex:multi3}\label{fig7}}
\end{figure}

Interaction of multi-fronts are shown in Fig.~\ref{fig7} for the set of parameters $(\alpha, \beta, \gamma, d, \tau, \theta)=(100, 100, -0.25, 5, 1, 1)$, and for the initial conditions
\[
u(x,0) = \left\{ \begin{array}{rcl}
         0 &,& \mbox{if }  x \in [50, 50],\\
          -1 &,& \mbox{otherwise}.\end{array} \right. \; , \quad v(x,0)=s(x,0)=0.
\]

\section{Pattern formation for the FitzHugh-Nagumo equation}
\label{sec:turing}

In this section, we will study Turing's pattern formation of FHN equation \eqref{fhn}. In order to determine what kind of spatial patterns is formed for Turing systems, an analytic criterion is developed in \citep{lago14}. Here, we will give this result and Turing instability conditions for the FHN  equation \eqref{fhn} with two different sets of $f_1(u)$, $f_2(v)$, $g_1(u)$, $g_2(v)$ in two-dimensional spatial domain $ \Omega$.
Alan Turing proposed this relationship between the parameters of reaction-diffusion equations and pattern formation \citep{Turing52}. He showed that under certain conditions any stable steady state can be converted into an unstable state. This is known as \textit{diffusion driven instability} or \textit{Turing instability}. The Turing instability conditions for the system \eqref{fhn} are given by
\begin{equation} \label{turing}
\begin{aligned}
\frac{\partial f_1}{\partial u} + \frac{\partial g_2}{\partial v} &< 0,  \\
\frac{\partial f_1}{\partial u} \cdot \frac{\partial g_2}{\partial v} - \frac{\partial f_2}{\partial v} \cdot \frac{\partial g_1}{\partial u}   &> 0,  \\
d_2 \frac{\partial f_1}{\partial u} + d_1 \frac{\partial g_2}{\partial v} &> 0,  \\
\left(d_2 \frac{\partial f_1}{\partial u} + d_1 \frac{\partial g_2}{\partial v}\right)^2 &> 4 d_1 d_2 \left(\frac{\partial f_1}{\partial u} \cdot \frac{\partial g_2}{\partial v} - \frac{\partial f_2}{\partial v} \cdot \frac{\partial g_1}{\partial u}\right),
\end{aligned}
\end{equation}
which have to be evaluated at the steady state solutions. The first two conditions of \eqref{turing} are already the stability conditions of the steady state solutions of \eqref{fhn}. The last two conditions make this stable steady state solutions unstable under the presence of the diffusion terms.

Turing patterns for FHN equation occur in form of spots or labyrinth-like patterns, depending on the number of global minimums of the Lyapunov energy functional \eqref{lyp} (see \citep{lago14}). When the  Lyapunov energy functional \eqref{lyp} of the system has a unique global minimum, then patterns appear as spots, otherwise, labyrinth-like patterns occur.. 

\subsection{Pattern examples}
\label{ex:turing}

We  fix the parameters  $\tau_1 = \tau_2 =1$, $d_1 = 0.00028$, and we take  the bi-stable nonlinear term $f_1(u)= u-u^3$, $f_2(v)=\kappa-v,$ $g_1(u)=u$ and $g_2(v)=-v$. First of all, note that the reaction terms include an even term $\kappa$. So for $\kappa = 0$, the pattern formation will be labyrinth-like patterns, and for nonzero $\kappa$ the results will be spots. We have to determine the values of $d_2$ so that the  Turing instability conditions \eqref{turing} are satisfied. The only steady state solution of \eqref{fhn} with the given conditions in this problem is $(u_0, v_0) \approx (-0.368403, -0.368403)$. When we evaluate the first order partial derivatives at this point, we get
$$
\frac{\partial f_1}{\partial u} \approx 0.592838,\quad \frac{\partial f_2}{\partial v} = -1,\quad \frac{\partial g_1}{\partial u} = 1,\quad \frac{\partial g_2}{\partial v} = -1.
$$
The first two conditions of \eqref{turing} are clearly satisfied. 
 In order to satisfy the third and last  conditions of \eqref{turing}, we need $d_2 > 0.000472$ and $d_2 > 0.002242$, respectively. Hence, for Turing instability, we obtain $d_2> 0.002242$ .

\begin{figure}[htb!]
\centering
\subfloat[]{\includegraphics[scale=0.25]{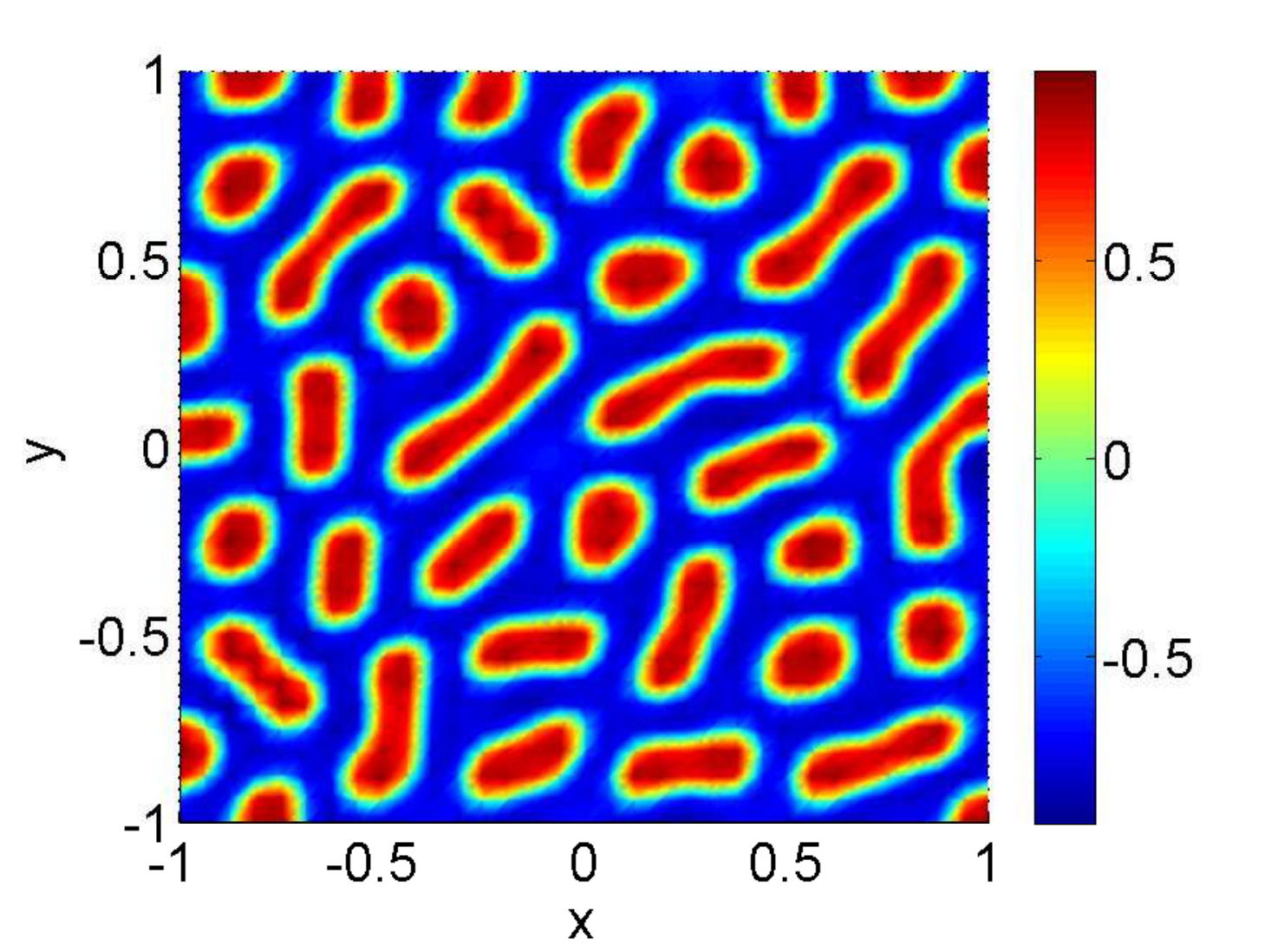}}
\subfloat[]{\includegraphics[scale=0.25]{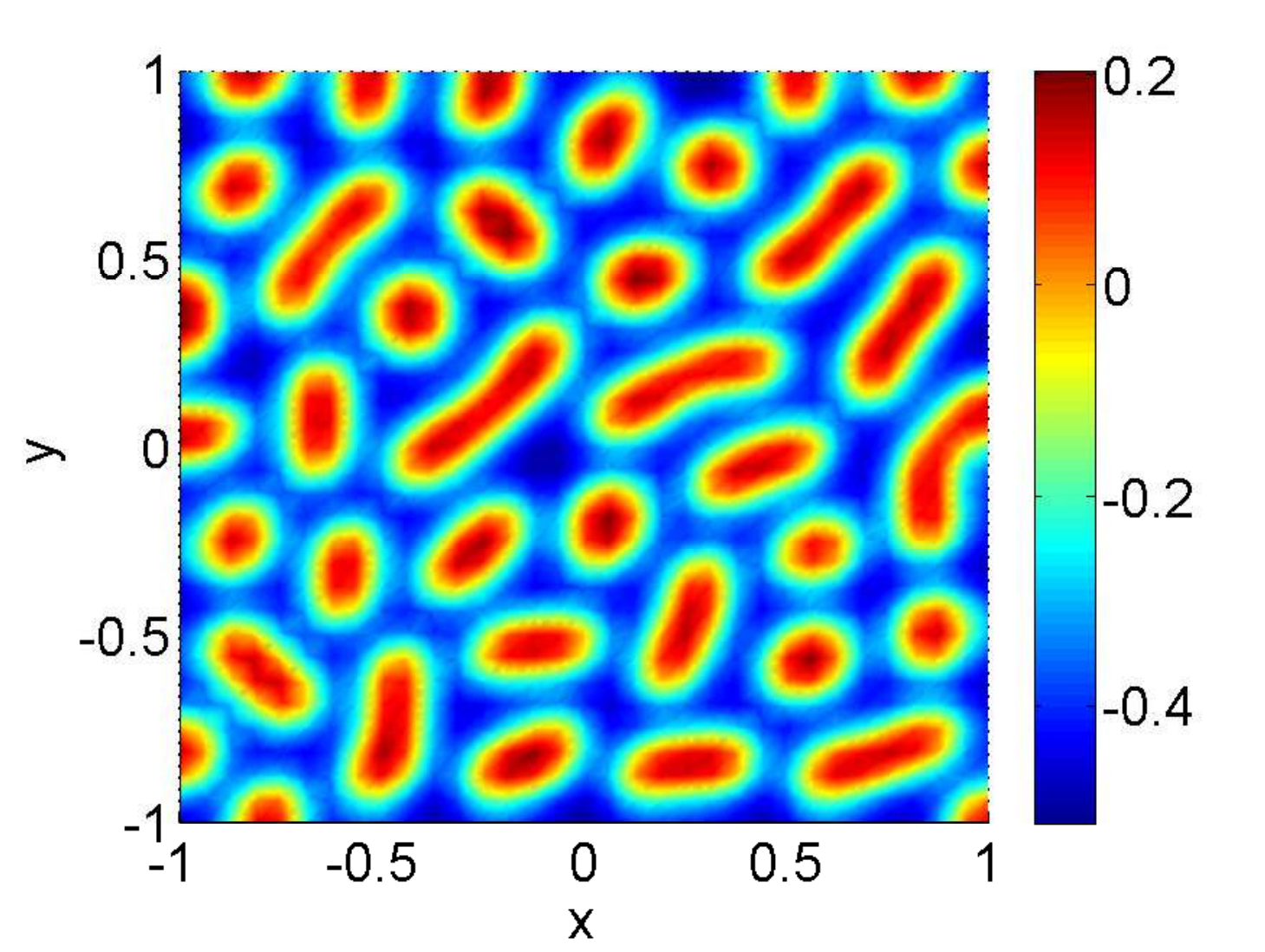}}
\subfloat[]{\includegraphics[scale=0.25]{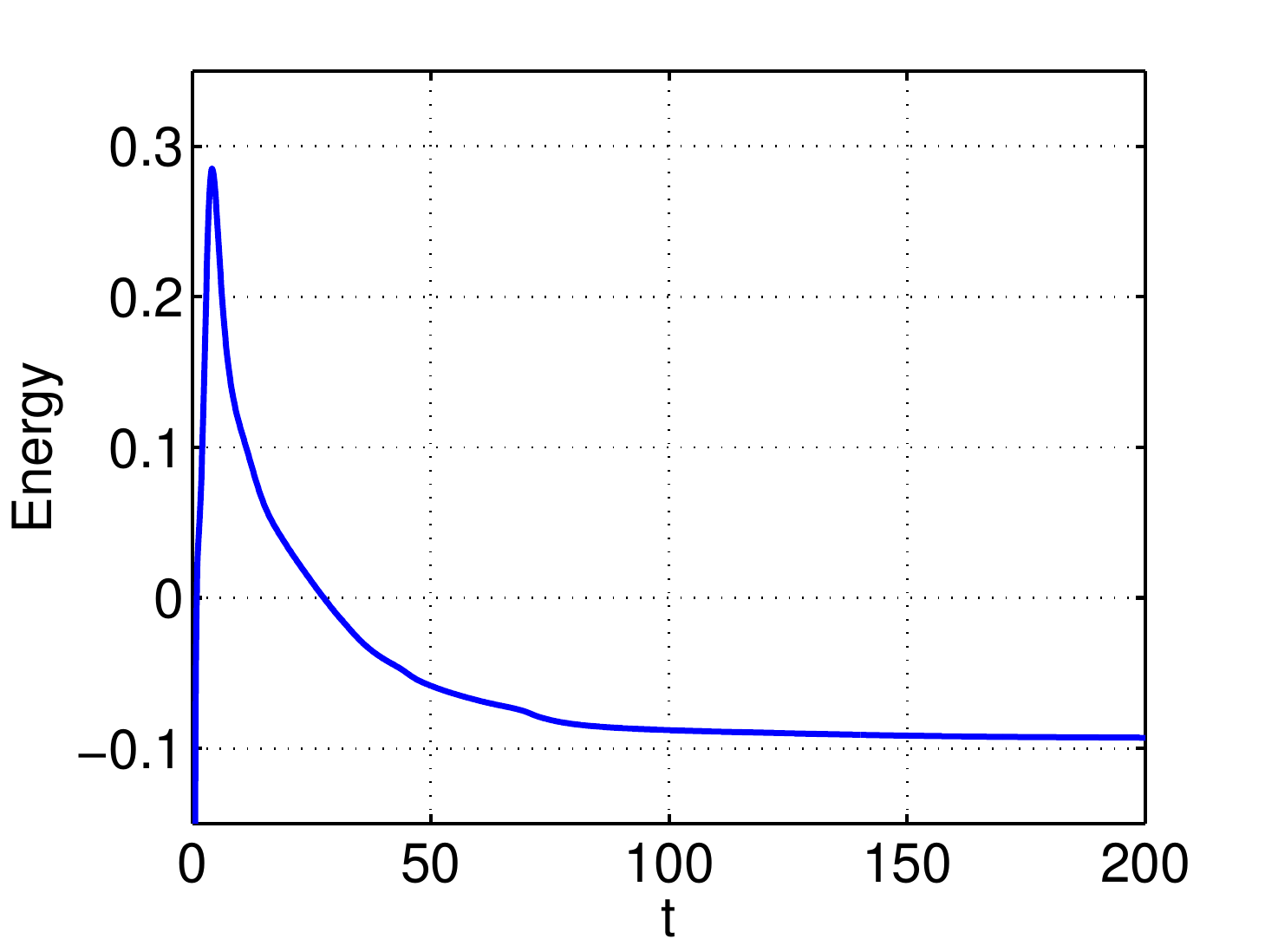}}
\caption{Spots: Solution profile of the first component u (left), second component v (middle) and evolution of discrete energy (right) for Example \ref{ex:turing}}
\label{pfig6}
\end{figure}

\begin{figure}[htb!]
\centering
\subfloat[]{\includegraphics[scale=0.25]{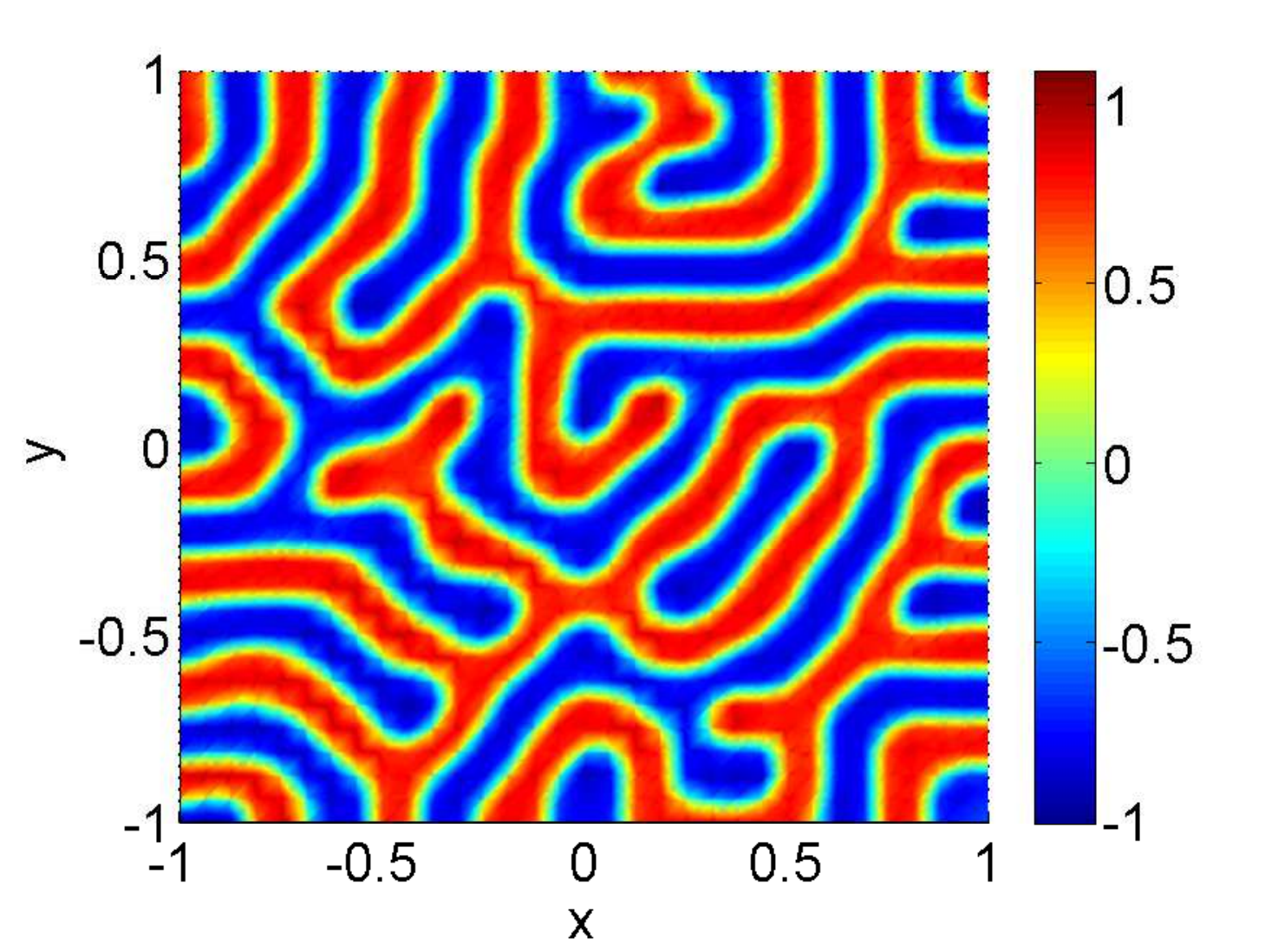}}
\subfloat[]{\includegraphics[scale=0.25]{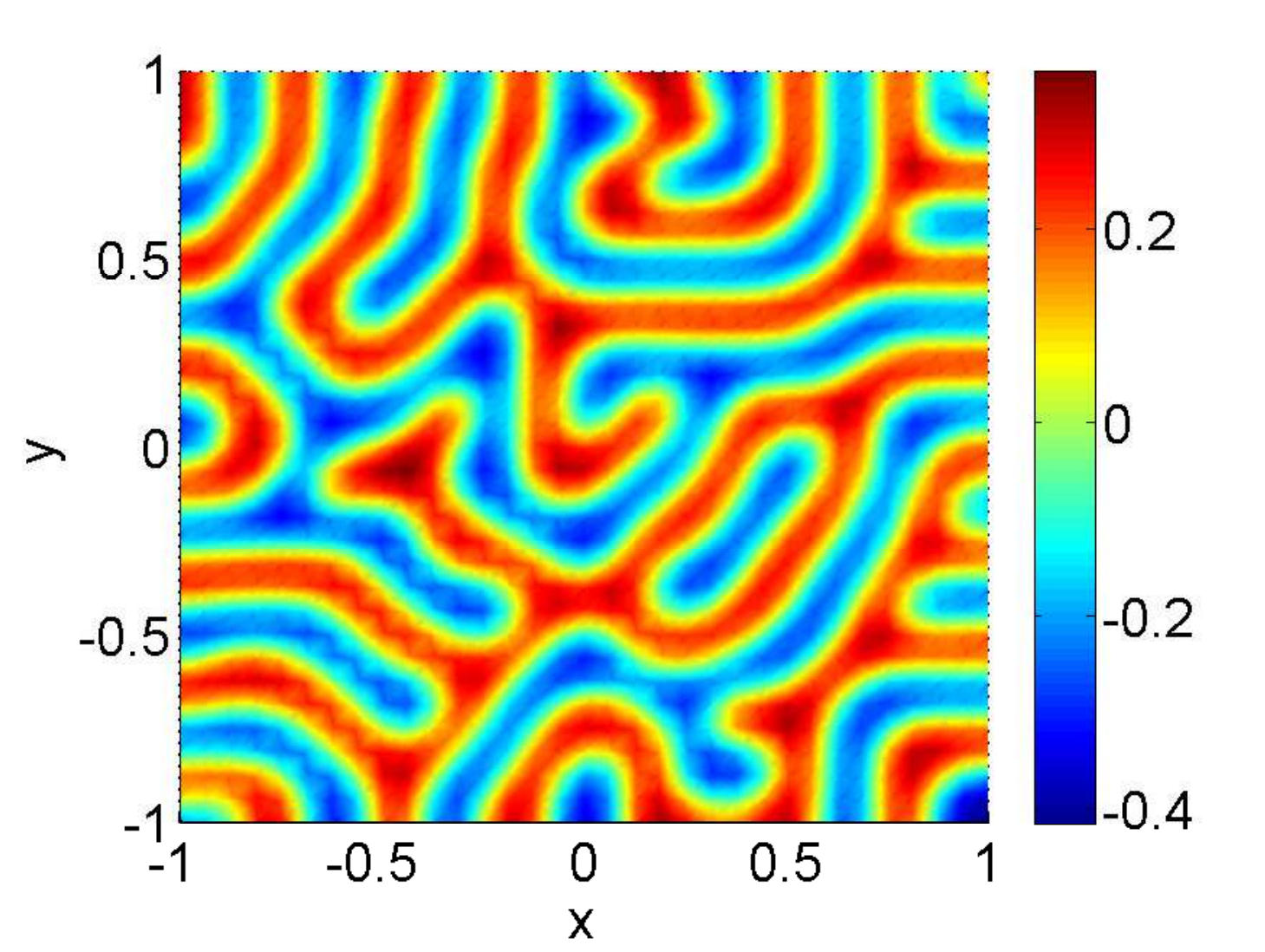}}
\subfloat[]{\includegraphics[scale=0.25]{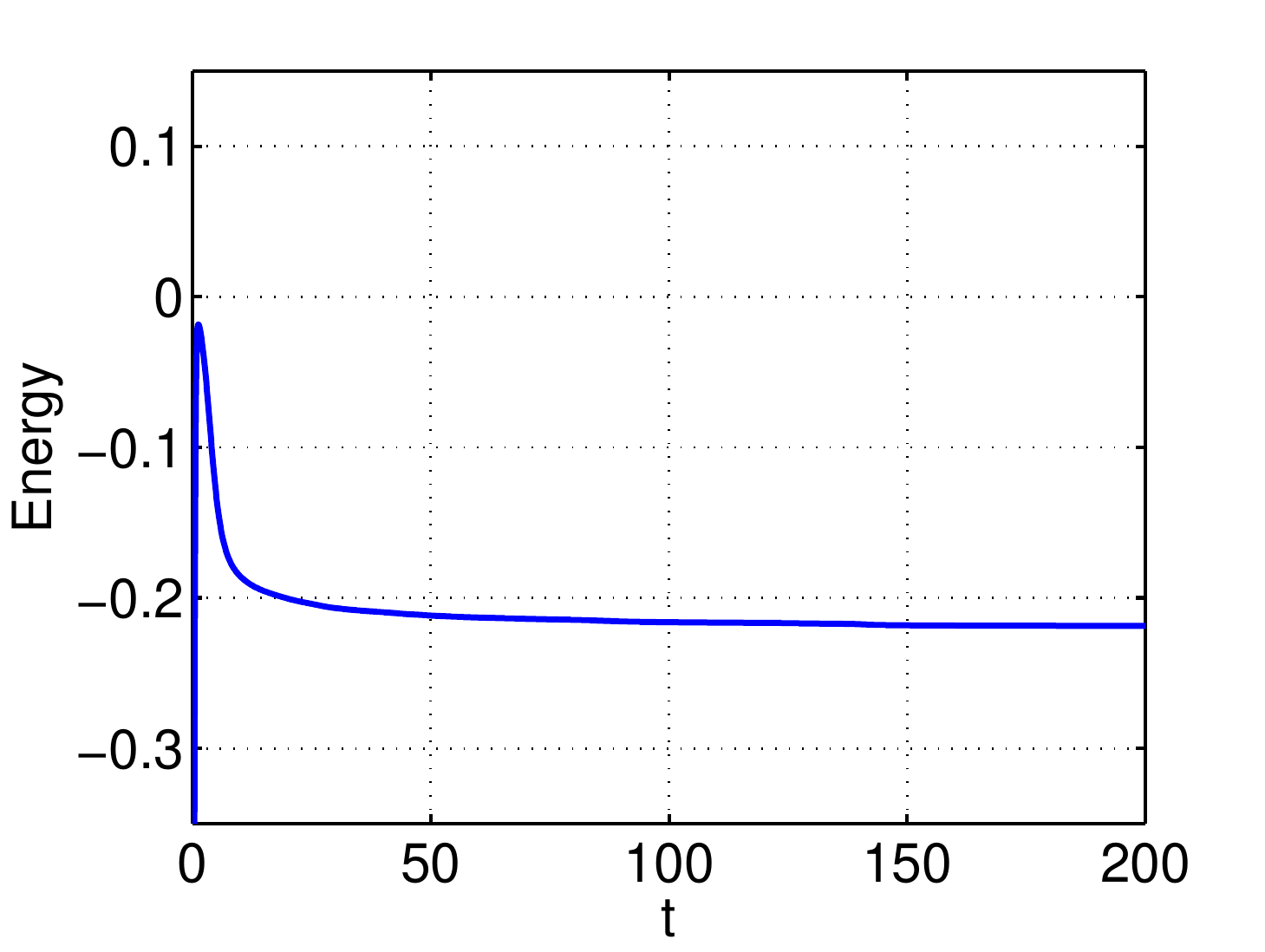}}
\caption{Labyrinth-like patterns: Solution profile of the first component u (left), second component v (middle) and evolution of discrete energy (right)  for Example \ref{ex:turing}}
\label{fig9}
\end{figure}

We take $d_2= 0.005$ and consider the domain $(x,y) \in [-1,1] \times [-1,1]$ and $t \in [0,200]$ with the uniformly distributed random number between $-1$ and $1$. Space mesh size and time step size are taken as $\Delta x=1/16$ and $\Delta t =1/10$, respectively.  The steady state is reached at $t=200$. The corresponding solution profiles and energy plots of spot and labyrinth-like patterns are given in Fig.~\ref{pfig6} and Fig.~\ref{fig9}, respectively.

\section{Reduced order model approximation for the FitzHugh-Nagumo system}
\label{sec:rom}

For MOR we consider the following special form of the FHN equation \eqref{fhn}
\begin{equation}\label{fhn_cont}
\begin{aligned}
u_t  &=  d_1\Delta u + f(u) - v + \kappa\\
v_t  &=  d_2\Delta v - v + u
\end{aligned}
\end{equation}
The SIPG semi-discretization of(\ref{fhn_cont}) leads to the following system of ODEs
\begin{equation}\label{fhn_ode}
\begin{aligned}
M\mathbf{u}_t  &=  S_u\mathbf{u} + F(\mathbf{u}) - M\mathbf{v} + K\\
M\mathbf{v}_t  &=  S_v\mathbf{v} - M\mathbf{v} + M\mathbf{u}
\end{aligned}
\end{equation}
with the stiffness matrices $S_u,S_v\in\mathbb{R}^{N\times N}$, $M\in\mathbb{R}^{N\times N}$ is the mass matrix, $K\in\mathbb{R}^{N}$ is the constant vector related to the integral containing the parameter $\kappa$ and $F\in\mathbb{R}^{N}$ is the vector of the bi-stable nonlinear terms of $\mathbf{u}=(u_1,u_2,\ldots , u_{N})$. The number $N=N_{loc}\times N_{el}$ stands for the degrees of freedom in DG methods, where $N_{loc}$ and $N_{el}$ represent the number of local dimension in each element (triangle) and the number of elements, respectively. The solutions of \eqref{fhn_ode} are of the form
\begin{equation}\label{fhn_comb}
u(t)=\sum_{i=1}^{N}u_i(t)\phi_i(x)=\phi\mathbf{u}(t) \; , \quad v(t)=\sum_{i=1}^{N}v_i(t)\phi_i(x)=\phi\mathbf{v}(t)
\end{equation}
where $u_i(t),v_i(t)\in\mathbb{R}$ are the unknown coefficients and $\phi_i$ are the DG basis functions.

Because the computation of the Turing patterns for \eqref{fhn_ode} is very time consuming,  we construct the reduced order model (ROM) of small dimension \eqref{fhn_ode} by utilizing the POD method \citep{Kunisch01}.


\subsection{Reduced order model}

The ROM  for  \eqref{fhn_ode} of dimension $k\ll N$ is formed by approximating the solutions $u(t)$ and $v(t)$ of the FOM \eqref{fhn_ode} in a  subspaces spanned by a set of $M$-orthogonal  basis functions $\{\psi_{u,i}\}_{i=1}^k$ and $\{\psi_{v,i}\}_{i=1}^k$ of dimension $k$ in $\mathbb{R}^{N}$ and then projecting onto that subspace. The approximate ROM solutions have the form
\begin{equation}\label{fhn_combrom}
u(t)\approx \sum_{i=1}^{k} \tilde{u}_i(t) \psi_{u,i} \; , \quad v(t)\approx \sum_{i=1}^{k} \tilde{v}_i(t) \psi_{v,i},
\end{equation}
where $\mathbf{\tilde{u}}(t)=(\tilde{u}_1(t),\ldots ,\tilde{u}_k(t))^T$ and $\mathbf{\tilde{v}}(t)=(\tilde{v}_1(t),\ldots ,\tilde{v}_k(t))^T$ are the coefficient vectors of the ROM solutions. The reduced basis functions $\{\psi_{u,i}\}$ and $\{\psi_{v,i}\}$ of the ROM solutions in \eqref{fhn_combrom} are constructed with DG finite element basis functions $\{\phi_i\}_{i=1}^{N}$
\begin{equation}\label{fhn_combpod}
\psi_{u,i} = \sum_{j=1}^{N} \Psi_{u,j,i} \phi_j(x)=\phi\Psi_u\; , \quad \psi_{v,i} = \sum_{j=1}^{N} \Psi_{v,j,i} \phi_j(x)=\phi\Psi_v,
\end{equation}
where the coefficient vectors of the $i-th$ reduced basis functions $\psi_{u,i}$ and $\psi_{v,i}$ are located at the $i-th$ columns of the matrices $\Psi_u=[\Psi_{u,\cdot ,1}, \ldots ,\Psi_{u,\cdot ,k} ]\in\mathbb{R}^{N\times k}$ and $\Psi_v=[\Psi_{v,\cdot ,1}, \ldots ,\Psi_{v,\cdot ,k} ]\in\mathbb{R}^{N\times k}$, respectively.

The M-orthogonal reduced basis functions $\{\psi_{u,i}\}$ and $\{\psi_{v,i}\}$, $i=1,2,\ldots , k$, are computed by the POD \citep{Kunisch01}. We consider the snapshot matrices $\mathcal{U}=[\mathbf{u}^1,\ldots ,\mathbf{u}^J]$ and $\mathcal{V}=[\mathbf{v}^1,\ldots ,\mathbf{v}^J]$ in $\mathbb{R}^{N\times J}$, where each member of the snapshot matrices $\mathcal{U}$ and $\mathcal{V}$ are the corresponding coefficient vectors of the discrete solutions $\{ u^i\}_{i=1}^J$ and $\{ v^i\}_{i=1}^J$, respectively, of the FOM \eqref{fhn_ode} at the time  $t_i$, $i=0,1,\ldots , J$, with $u^i\approx u(t_i)$ and $v^i\approx v(t_i)$.
Then, for $w\in\{u,v\}$, the M-orthogonal reduced basis functions $\{\psi_{w,i}\}$, $i=1,2,\ldots , k$, are given by the solution of the minimization problem
\begin{align*}
\min_{\psi_{w,1},\ldots ,\psi_{w,k}} \frac{1}{J}\sum_{j=1}^J \left\| w^j - \sum_{i=1}^k (w^j,\psi_{w,i})_{L^2(\Omega)}\psi_{w,i}\right\|_{L^2(\Omega)}^2 \\
\text{subject to } (\psi_{w,i},\psi_{w,j})_{L^2(\Omega)} = \Psi_{w,\cdot ,i}^TM\Psi_{w,\cdot ,j}=\delta_{ij} \; , \; 1\leq i,j\leq k,
\end{align*}
where $\delta_{ij}$ is the Kronecker delta. The above minimization problem is equivalent to the eigenvalue problem
\begin{equation}\label{eg1}
\mathcal{U}\mathcal{U}^TM\Psi_{u,\cdot ,i}=\sigma_{u,i}^2\Psi_{u,\cdot ,i} \; , \qquad \mathcal{V}\mathcal{V}^TM\Psi_{v,\cdot ,i}=\sigma_{v,i}^2\Psi_{v,\cdot ,i} \; , \quad i=1,2,\ldots ,k
\end{equation}
for the coefficient vectors $\Psi_{u,\cdot ,i}$ and $\Psi_{v,\cdot ,i}$ of the POD basis functions $\psi_{u,i}$ and $\psi_{v,i}$, respectively. Defining $\widehat{\mathcal{U}}=R\mathcal{U}$ and $\widehat{\mathcal{V}}=R\mathcal{V}$ ($R^T$ is the Cholesky factor of the mass matrix $M$), we obtain the equivalent formulation of  \eqref{eg1} as
\begin{equation}\label{eg2}
\widehat{\mathcal{U}}\widehat{\mathcal{U}}^T\widehat{\Psi}_{u,\cdot ,i}=\sigma_{u,i}^2\widehat{\Psi}_{u,\cdot ,i} \; , \qquad \widehat{\mathcal{V}}\widehat{\mathcal{V}}^T\widehat{\Psi}_{v,\cdot ,i}=\sigma_{v,i}^2\widehat{\Psi}_{v,\cdot ,i} \; , \quad i=1,2,\ldots ,k
\end{equation}
where $\widehat{\Psi}_{\cdot,\cdot ,i}=R\Psi_{\cdot,\cdot ,i}$. The solutions $\widehat{\Psi}_{\cdot,\cdot ,i}$  of \eqref{eg2} are obtained as the first $k$ left singular vectors $\widehat{\Psi}_{u,\cdot ,i}=\zeta_{u,i}$ and $\widehat{\Psi}_{v,\cdot ,i}=\zeta_{v,i}$ in the singular value decomposition (SVD) of $\widehat{\mathcal{U}}$ and $\widehat{\mathcal{V}}$, respectively, as
$$
\widehat{\mathcal{U}} = \zeta_{u}\Sigma_u \beta_u^T \; , \qquad \widehat{\mathcal{V}} = \zeta_{v}\Sigma_v \beta_v^T,
$$
where the diagonal matrices $\Sigma_u$ and $\Sigma_v$ contain the singular values $\sigma_{u,i}$ and $\sigma_{v,i}$ on the diagonals, respectively. Using  $\widehat{\Psi}_{\cdot,\cdot ,i}=R\Psi_{\cdot,\cdot ,i}$, the coefficient vectors $\Psi_{\cdot,\cdot ,i}$ of the POD basis functions are computed as
$$
\Psi_{u,\cdot ,i}=R^{-1}\widehat{\Psi}_{u,\cdot ,i} \; , \quad \Psi_{v,\cdot ,i}=R^{-1}\widehat{\Psi}_{v,\cdot ,i} \; , \quad i=1,2,\ldots k.
$$
In addition, using the expansions \eqref{fhn_comb}, \eqref{fhn_combrom} and \eqref{fhn_combpod}, we have
\begin{align}
\mathbf{u} = \Psi_u\mathbf{\tilde{u}} \; , \quad \mathbf{v} = \Psi_v\mathbf{\tilde{v}}.\label{fulltorom}
\end{align}

For the construction of the $k$-dimensional ROM, we substitute  the relations \eqref{fulltorom} into the system \eqref{fhn_ode} and we project onto the reduced spaces spanned by $\{\psi_{u,1},\ldots ,\psi_{u,k}\}$  and $\{\psi_{v,1},\ldots ,\psi_{v,k}\}$, respectively, leading to the system
\begin{equation}\label{fhn_rom}
\begin{aligned}
\mathbf{\tilde{u}}_t  &=  \tilde{A}\mathbf{\tilde{u}} + \Psi_u^TF(\Psi_u\mathbf{\tilde{u}}) - \tilde{M}_u\mathbf{\tilde{v}} + \Psi_u^TK\\
\mathbf{\tilde{v}}_t  &=  \tilde{B}\mathbf{\tilde{v}} - \mathbf{\tilde{v}} + \tilde{M}_v\mathbf{\tilde{u}}
\end{aligned}
\end{equation}
with the reduced matrices
$$
\tilde{A}=\Psi_u^TA\Psi_u \; , \quad \tilde{B}=\Psi_v^TB\Psi_v \; , \quad \tilde{M}_u=\Psi_u^TM\Psi_v\; , \quad \tilde{M}_v=\Psi_v^TM\Psi_u.
$$



\subsection{Discrete empirical interpolation method (DEIM)}

Although the dimension of the reduced system \eqref{fhn_rom} is small, $k\ll N$, the computation of the nonlinear term
\begin{equation}\label{podN}
N(\mathbf{\tilde{u}})=\Psi_u^TF(\Psi_u\mathbf{\tilde{u}})
\end{equation}
still depends on the dimension $N$ of the full system. We apply the DEIM \citep{chaturantabut10nmr} to approximate the nonlinear function $F(\Psi_u\mathbf{\tilde{u}})$ from a subspace generated by the non-linear functions. Let $\mathcal{F}=[F_1,F_2,\ldots , F_J]\in\mathbb{R}^{N\times J}$ denotes the snapshot matrix of the nonlinear functions at each time instances $t_1,\ldots ,t_J$ computed as the solution of the full system \eqref{fhn_ode} with $F_i=F(\Psi_u\mathbf{\tilde{u}}(t_i))$, $i=1,2,\ldots ,J$. Using the SVD of the matrix $\mathcal{F}$, we can find $m\ll N$ orthogonal  basis functions $\{W_i\}_{i=1}^m$ spanning the $m$-dimensional subspace of  $\mathcal{F}$. Then, with $W=[W_1,\ldots ,W_m]\in\mathbb{R}^{N\times m}$, we can approximate the nonlinear function by
\begin{equation}\label{podF}
F(\Psi_u\mathbf{\tilde{u}}(t)) \approx Ws(t)
\end{equation}
with the corresponding coefficient vector $s(t)$. We note that the system \eqref{podF} is overdetermined. Thus, to compute the coefficient vector $s(t)$, we take $m$ distinguished rows from the system $Ws(t)$ through the projection using a permutation matrix $P=[e_{\wp_1},\ldots , e_{\wp_m}]\in\mathbb{R}^{N\times m}$ where $e_{\wp_i}$ is the $i$-th column of the identity matrix $I\in\mathbb{R}^{N\times N}$ and computed by Algorithm \ref{deimalg} \citep{chaturantabut10nmr}  so that $P^TW$ is non-singular.
\begin{algorithm}
\caption{DEIM Algorithm\label{deimalg}}
\textbf{Input:} POD basis functions $\{W_i\}_{i=1}^m$
 \\
\textbf{Output:} Index Vector $\wp=[\wp_1,\ldots ,\wp_m]^T$ and Permutation Matrix $P=[e_{\wp_1},\ldots , e_{\wp_m}]$
 \\
\begin{algorithmic}[1]
\STATE \textbf{INPUT: } \{$W_i\}_{i=1}^m \subset \mathbb{R}^{N}$
\STATE \textbf{OUTPUT: } $\vec \wp = [\wp_1,...,\wp_m]^T \in \mathbb{R}^m, \ P \in \mathbb{R}^{N \times m}$
\STATE $[|\rho |,\wp_1] = \max \{|W_1|\}$
\STATE $W = [W_1], P = [\mathbf{e}_{\wp_1}], \vec \wp = [\wp_1]$
\FOR{$i = 2$ to $m$}
\STATE Solve $(P^T W) \mathbf{c} = P^T W_i$ for $\mathbf{c}$
\STATE $\mathbf{r} = W_i - W \mathbf{c}$
\STATE $[|\rho|,\wp_i] = \max \{|\mathbf{r}|\}$
\STATE $W \leftarrow [W \ W_i],  P \leftarrow [P \ \mathbf{e}_{\wp_i}], \vec \wp \leftarrow \begin{bmatrix} \vec \wp \\ \wp_i \end{bmatrix}$
\ENDFOR
\end{algorithmic}
\label{algorithm:linear}
\end{algorithm}

The projection of \eqref{podF} leads to the system
\begin{equation}\label{prF}
s(t)= (P^TW)^{-1}P^T F(\Psi_u\mathbf{\tilde{u}}(t)).
\end{equation}
Then, using \eqref{podF} and \eqref{prF},  the non-linear term \eqref{podN} can be approximated as
\begin{equation}\label{prN}
N(\mathbf{\tilde{u}}(t)) \approx \tilde{N}(\mathbf{\tilde{u}}(t)) = Q\tilde{F} ,
\end{equation}
where the matrix $Q=\Psi_u^TW(P^TW)^{-1}\in\mathbb{R}^{k\times m}$ is precomputable and $\tilde{F}=P^T F(\Psi_u\mathbf{\tilde{u}}(t))\in\mathbb{R}^{m}$ is the $m$-dimensional non-linear vector which can be computed in an efficient way. 

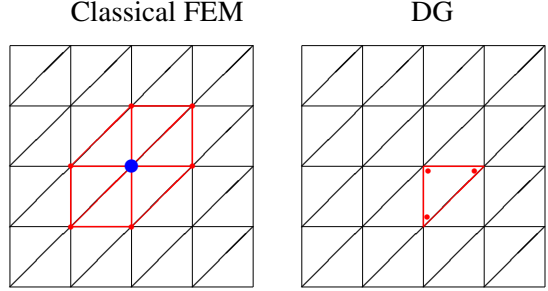
\begin{figure}[htb]
	\centering
	\setlength{\unitlength}{0.32mm}
		\begin{picture}(220, 120)
		
				\linethickness{0.1mm}
				\put(25,110){\large{Classical FEM}}
				\put(0,0){\line(1,0){100}}
				\put(0,0){\line(0,1){100}}
				\put(0,100){\line(1,0){100}}
				\put(100,0){\line(0,1){100}}
								
				\put(0,25){\line(1,0){100}}
				\put(0,50){\line(1,0){100}}
				\put(0,75){\line(1,0){100}}
				\put(25,0){\line(0,1){100}}
				\put(50,0){\line(0,1){100}}
				\put(75,0){\line(0,1){100}}
				
				\put(0,0){\line(1,1){100}}
				\put(0,25){\line(1,1){75}}
				\put(0,50){\line(1,1){50}}
				\put(0,75){\line(1,1){25}}
				\put(25,0){\line(1,1){75}}
				\put(50,0){\line(1,1){50}}
				\put(75,0){\line(1,1){25}}
				
				\color{red}
				\linethickness{0.2mm}
				\put(50,25){\line(1,1){25}}				
				\put(50,25){\line(0,1){25}}
				\put(50,50){\line(1,0){25}}			
				\put(25,25){\line(1,0){25}}	
				\put(25,25){\line(1,1){25}}	
				\put(25,25){\line(0,1){25}}	
				\put(25,50){\line(1,0){25}}	
				\put(25,50){\line(1,1){25}}	
				\put(50,50){\line(0,1){25}}		
				\put(50,50){\line(1,1){25}}		
				\put(50,75){\line(1,0){25}}		
				\put(75,50){\line(0,1){25}}	
					
				\put(25,25){\circle*{2.5}}
				\put(50,25){\circle*{2.5}}
				\put(25,50){\circle*{2.5}}
				\put(75,50){\circle*{2.5}}
				\put(50,75){\circle*{2.5}}
				\put(75,75){\circle*{2.5}}
				\color{blue}
				\put(50,50){\circle*{5}}

				\color{black}
				\linethickness{0.1mm}
				\put(165,110){\large{DG}}
				\put(120,0){\line(1,0){100}}
				\put(120,0){\line(0,1){100}}
				\put(120,100){\line(1,0){100}}
				\put(220,0){\line(0,1){100}}
				
				\put(120,25){\line(1,0){100}}
				\put(120,50){\line(1,0){100}}
				\put(120,75){\line(1,0){100}}
				\put(145,0){\line(0,1){100}}
				\put(170,0){\line(0,1){100}}
				\put(195,0){\line(0,1){100}}

				\put(120,0){\line(1,1){100}}
				\put(120,25){\line(1,1){75}}
				\put(120,50){\line(1,1){50}}
				\put(120,75){\line(1,1){25}}
				\put(145,0){\line(1,1){75}}
				\put(170,0){\line(1,1){50}}
				\put(195,0){\line(1,1){25}}
				
				\color{red}
				\linethickness{0.2mm}
				\put(170,25){\line(1,1){25}}				
				\put(170,25){\line(0,1){25}}
				\put(170,50){\line(1,0){25}}	
					
			  \put(171.5,29){\circle*{2.5}}
				\put(172,48){\circle*{2.5}}
				\put(191,48){\circle*{2.5}}
				
				\color{black}
				\linethickness{0.1mm}
		\end{picture}
		\caption{Connectivity of degrees of freedoms for linear basis\label{connect}}		
	\end{figure}

The DG requires only computation of the integrals on a single triangular element, which is not the case in continuous finite elements where all the interior degrees of freedoms are shared by usually 6 triangular elements (see Fig.~\ref{connect}).
The unassembled finite element approach is used in \citep{Heinkenschloss14}, so that each DEIM point is related to one element. This reduces the online computational cost,  but  increases the number of snapshots and therefore the cost of the offline computation.
Due to its' local nature, the DG discretization is automatically in the unassembled form and it does not require computation of additional snapshots (see Fig. \ref{connect}).

For the solution of the nonlinear system by  Newton's method, the entries of the Jacobian $J_F\in\mathbb{R}^{N\times N}$ of the non-linear terms are of the form  for $i,j=1,2,\ldots , N$
$$
(J_F)_{i,j}=\frac{\partial}{\partial u_j} F_i(\mathbf{u})=\int_{E_{t_i}}\frac{\partial f(\mathbf{u})}{\partial u_j}\phi_i dx = \int_{E_{t_i}}f'(\mathbf{u})\phi_j\phi_i dx.
$$

Because the DG basis functions are defined only a single element and they vanish outside that element, the integral terms of the Jacobian matrix $(J_F)_{i,j}$ vanish on the triangular elements $E_{t_i}$ where the basis function $\phi_j$ is not defined. Unlike the continuous finite elements where the Jacobian matrix contains overlapping blocks, Jacobian matrix in DG appears in block diagonal form, and has the form
$$
\frac{\partial}{\partial \mathbf{\tilde{u}}} N(\mathbf{\tilde{u}}) = \Psi_u^T J_F \Psi_u \; , \quad \frac{\partial}{\partial \mathbf{\tilde{u}}} \tilde{N}(\mathbf{\tilde{u}}) = Q (P^TJ_F) \Psi_u.
$$
We note that $(P^TJ_F)\in\mathbb{R}^{m\times N}$ is the matrix whose $i$-th row is the $\wp_i$-th row of the Jacobian $J_F$, $i=1,2,\ldots ,m$, and
in each row of the Jacobian there are only $N_{loc}$ nonzero terms because of the local structure of the DG.


\subsection{Numerical results }
We consider the FHN system \eqref{fhn_cont} with labyrinth-like patterns on $\Omega = [-1,1]^2$ with zero flux boundary conditions and with the initial conditions as uniformly distributed random numbers between $-1$ and $1$. We set the parameters $d_1=0.00028$, $d_2=0.005$ and $\kappa =0$.

\begin{figure}[htb!]
\centering
\includegraphics[scale=0.25]{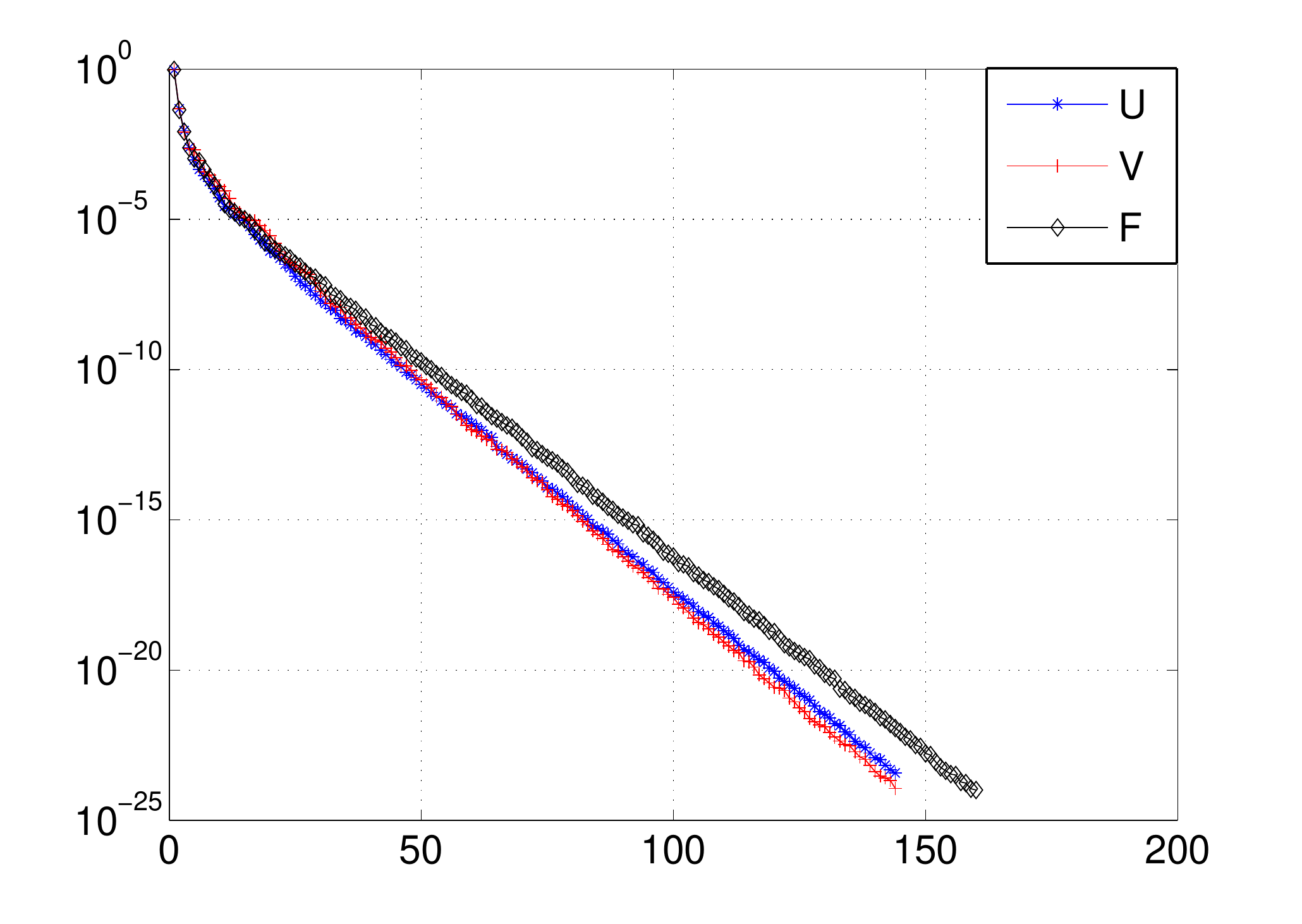}
\caption{Decay of the singular values related to $U$, $V$ and nonlinearity $F$\label{sec7_sing}}
\end{figure}

\begin{figure}[htb!]
\centering
\subfloat[]{\includegraphics[scale=0.25]{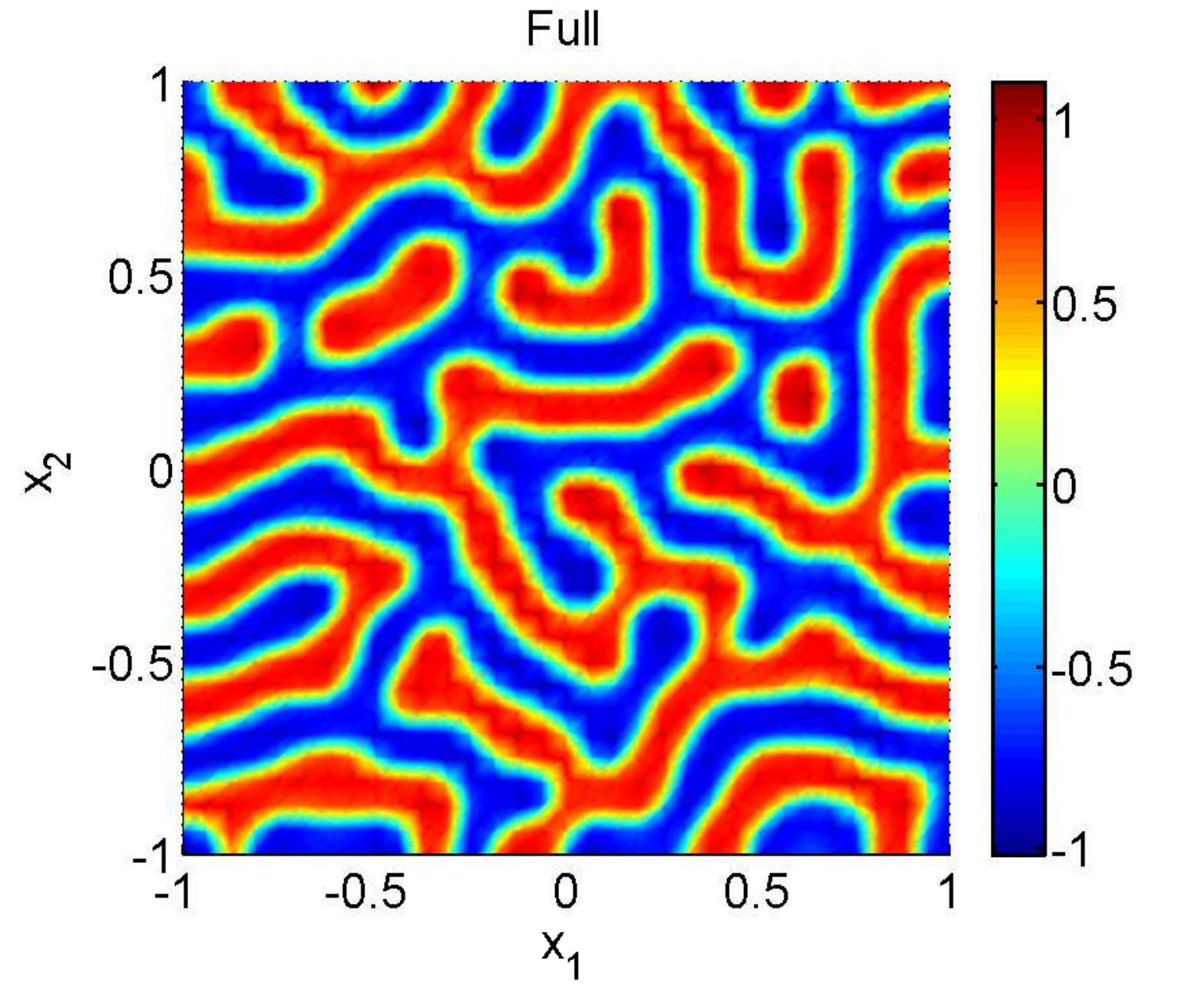}}
\subfloat[]{\includegraphics[scale=0.25]{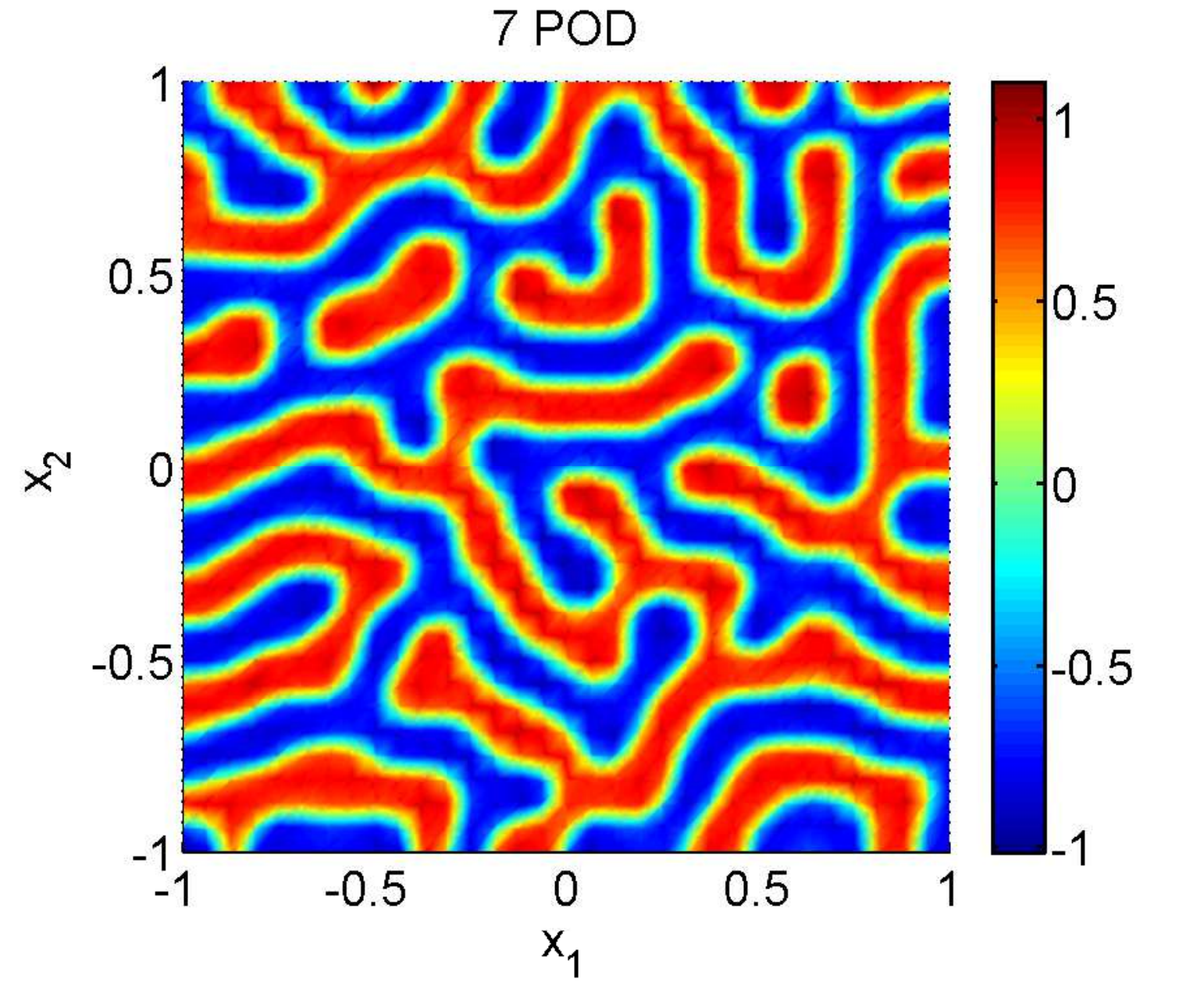}}
\subfloat[]{\includegraphics[scale=0.25]{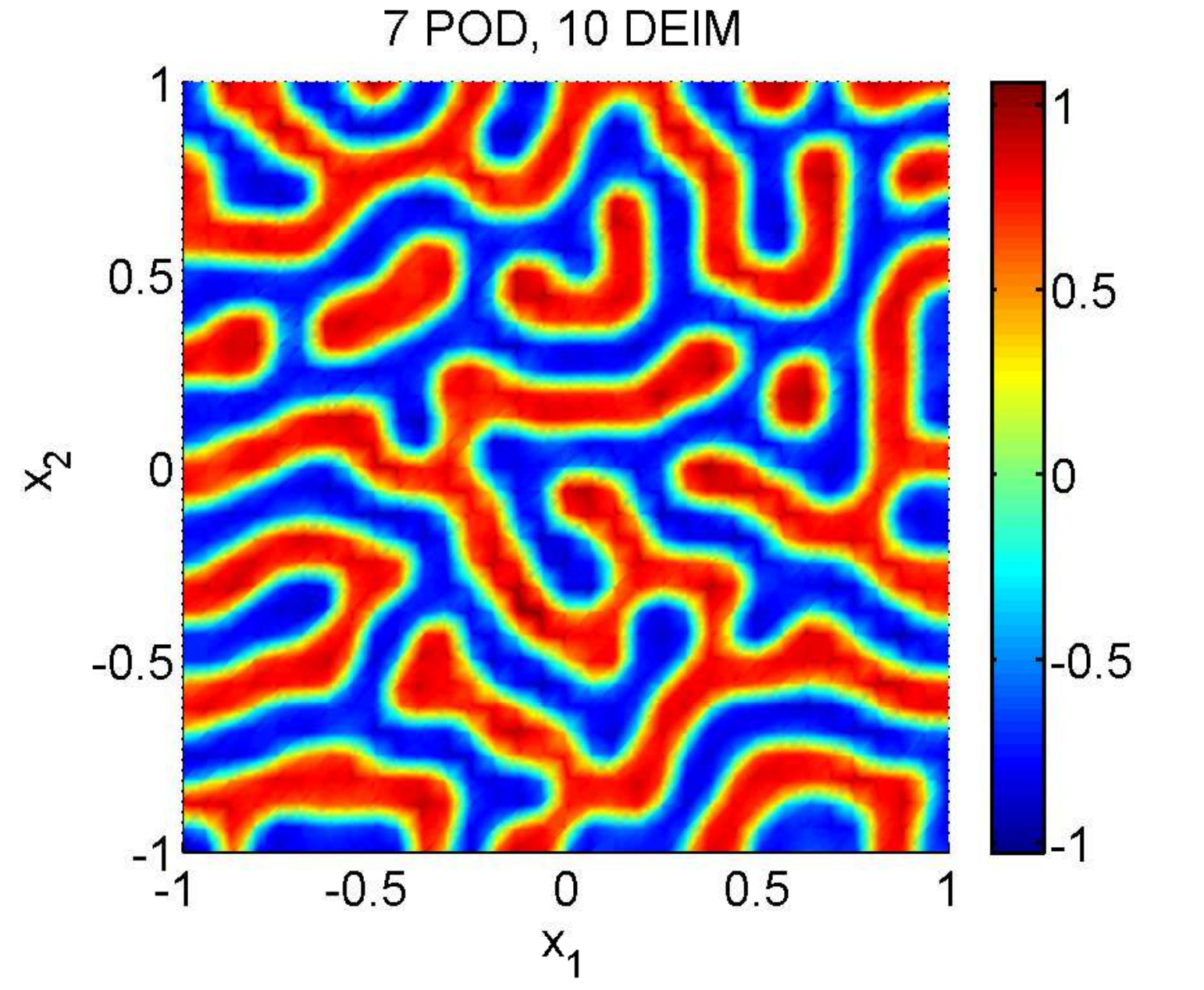}}
\caption{Solution profiles of the component $u$ at the steady state\label{sec7_plotu}}
\end{figure}

\begin{figure}[htb!]
\centering
\subfloat[]{\includegraphics[scale=0.25]{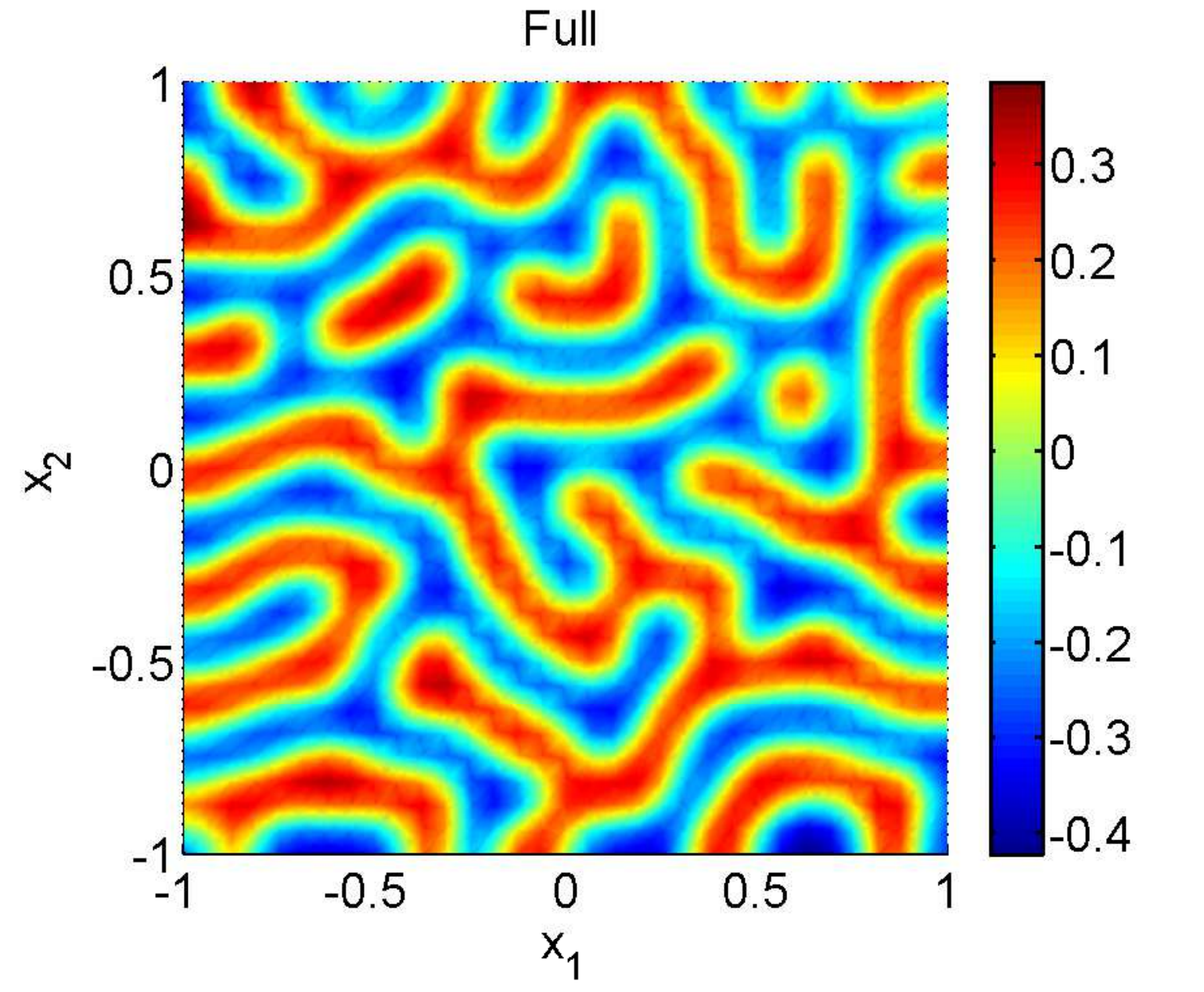}}
\subfloat[]{\includegraphics[scale=0.25]{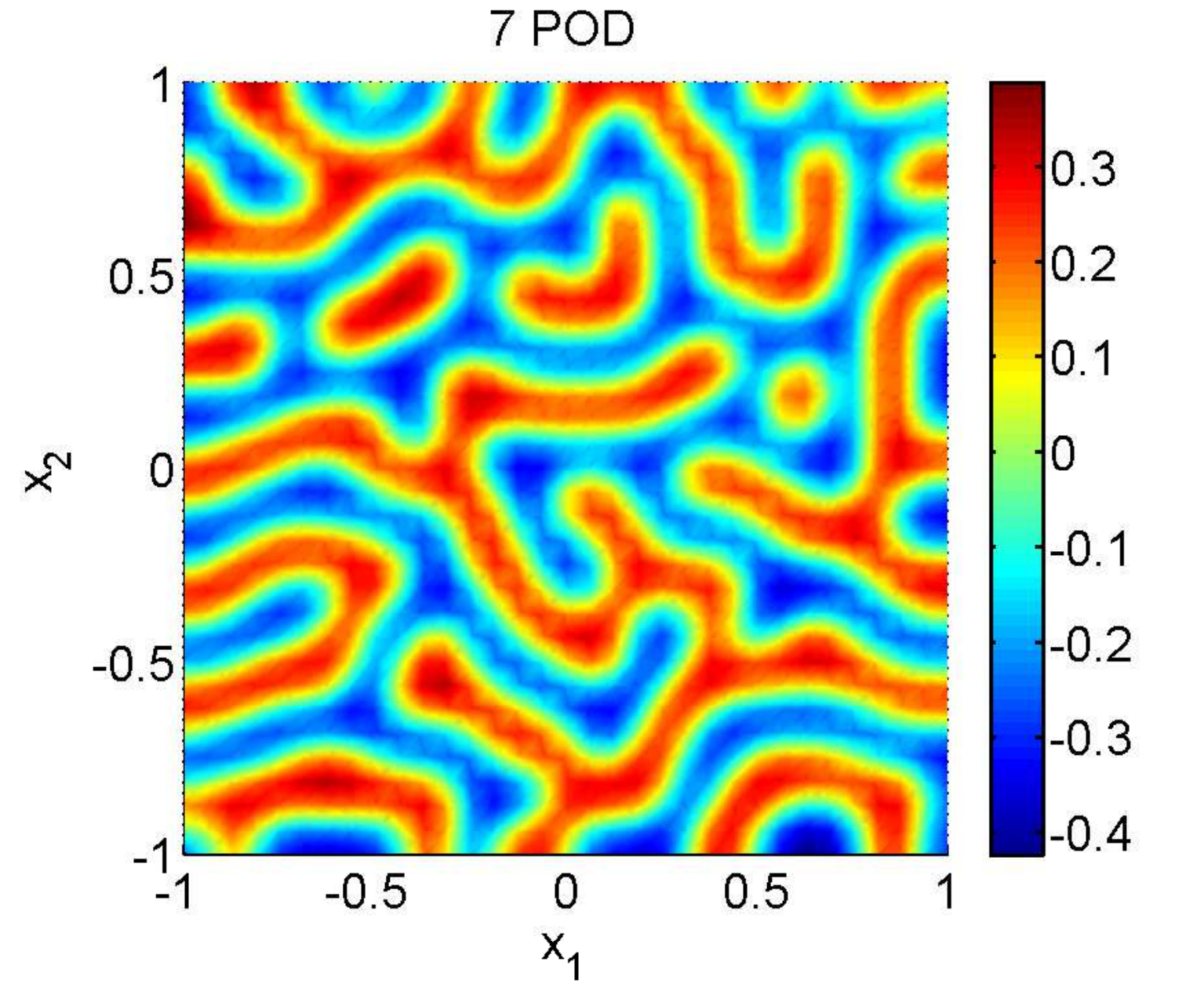}}
\subfloat[]{\includegraphics[scale=0.25]{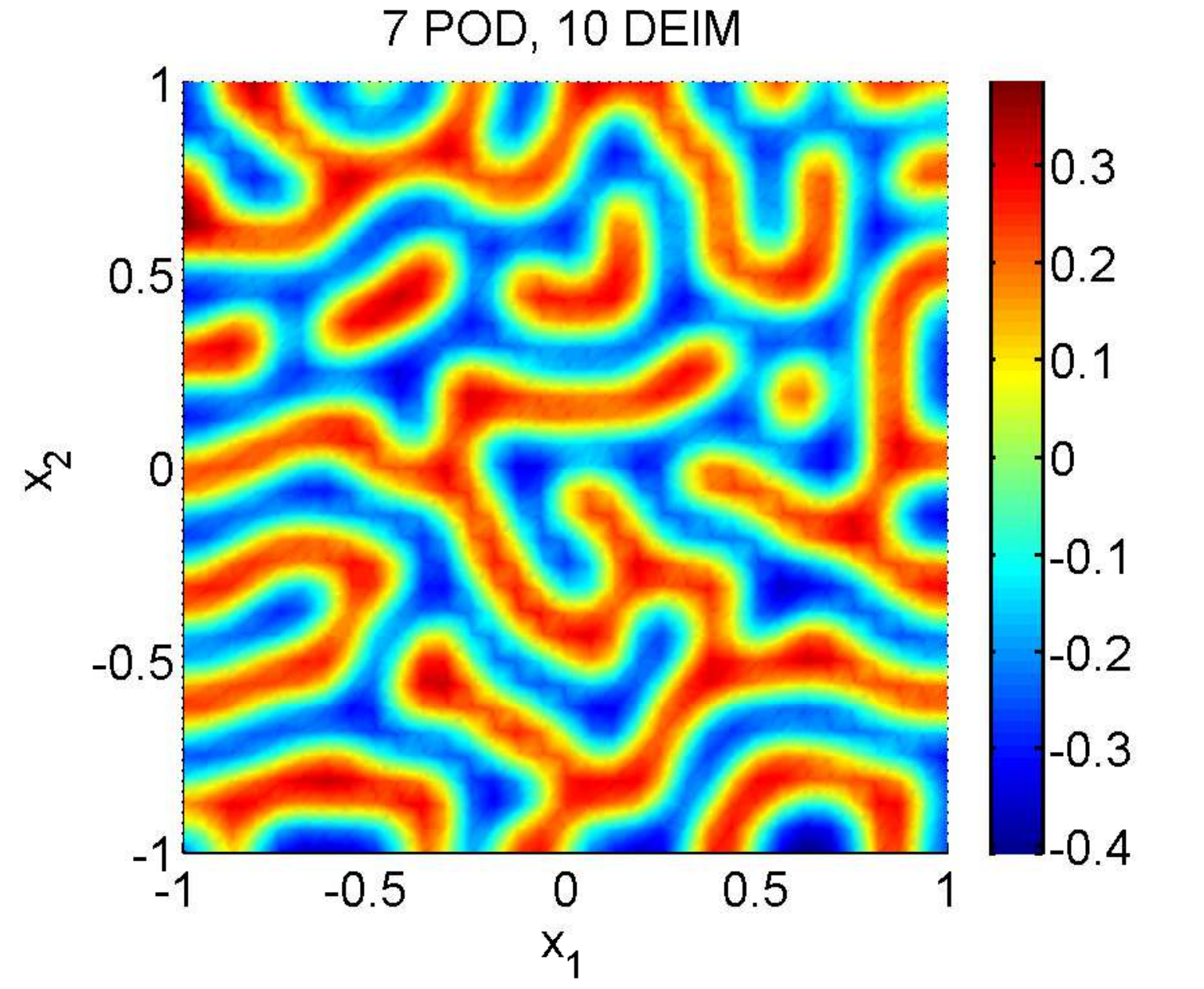}}
\caption{Solution profiles of the component $v$ at the steady state\label{sec7_plotv}}
\end{figure}

\begin{figure}[htb!]
\centering
\subfloat[]{\includegraphics[scale=0.25]{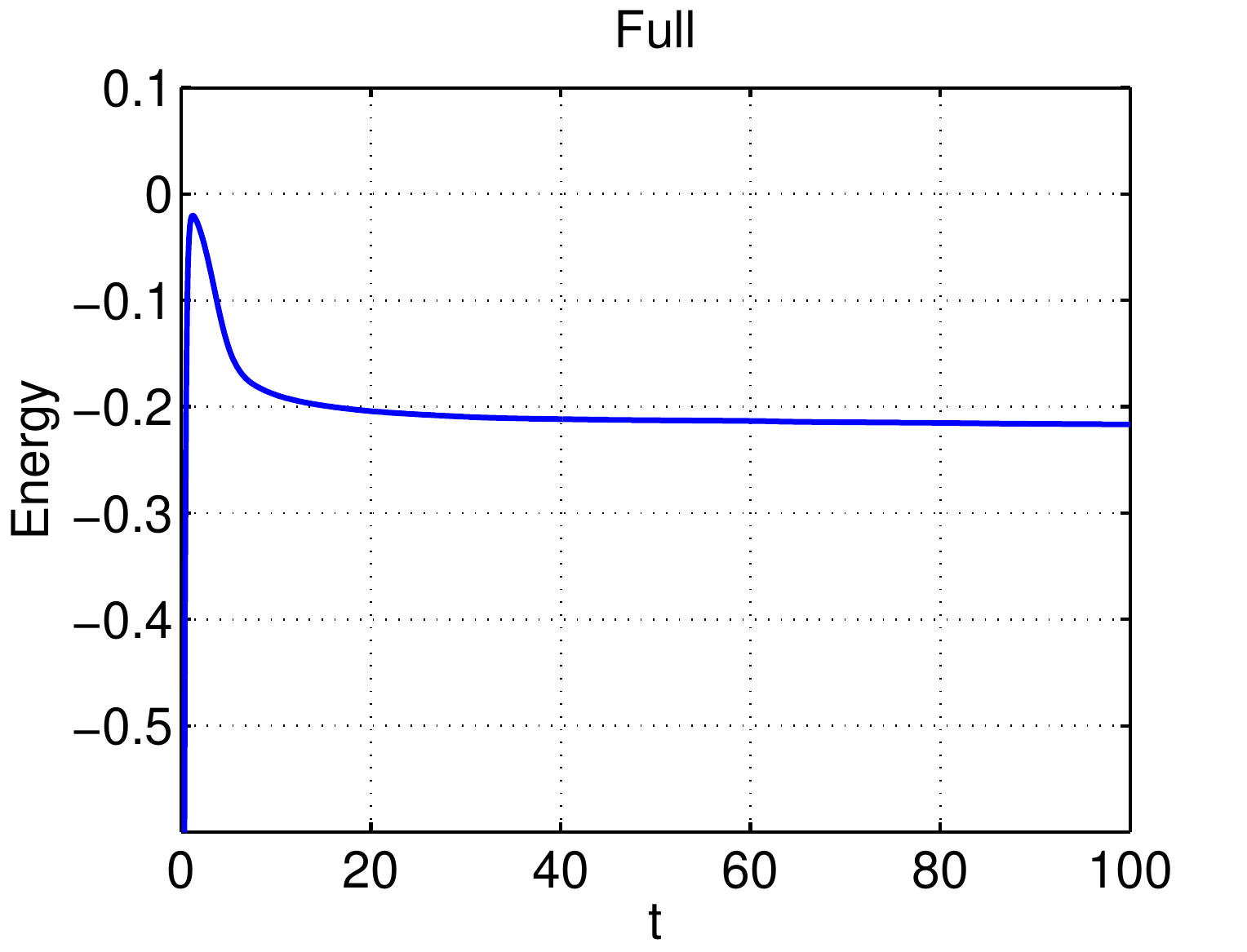}}
\subfloat[]{\includegraphics[scale=0.25]{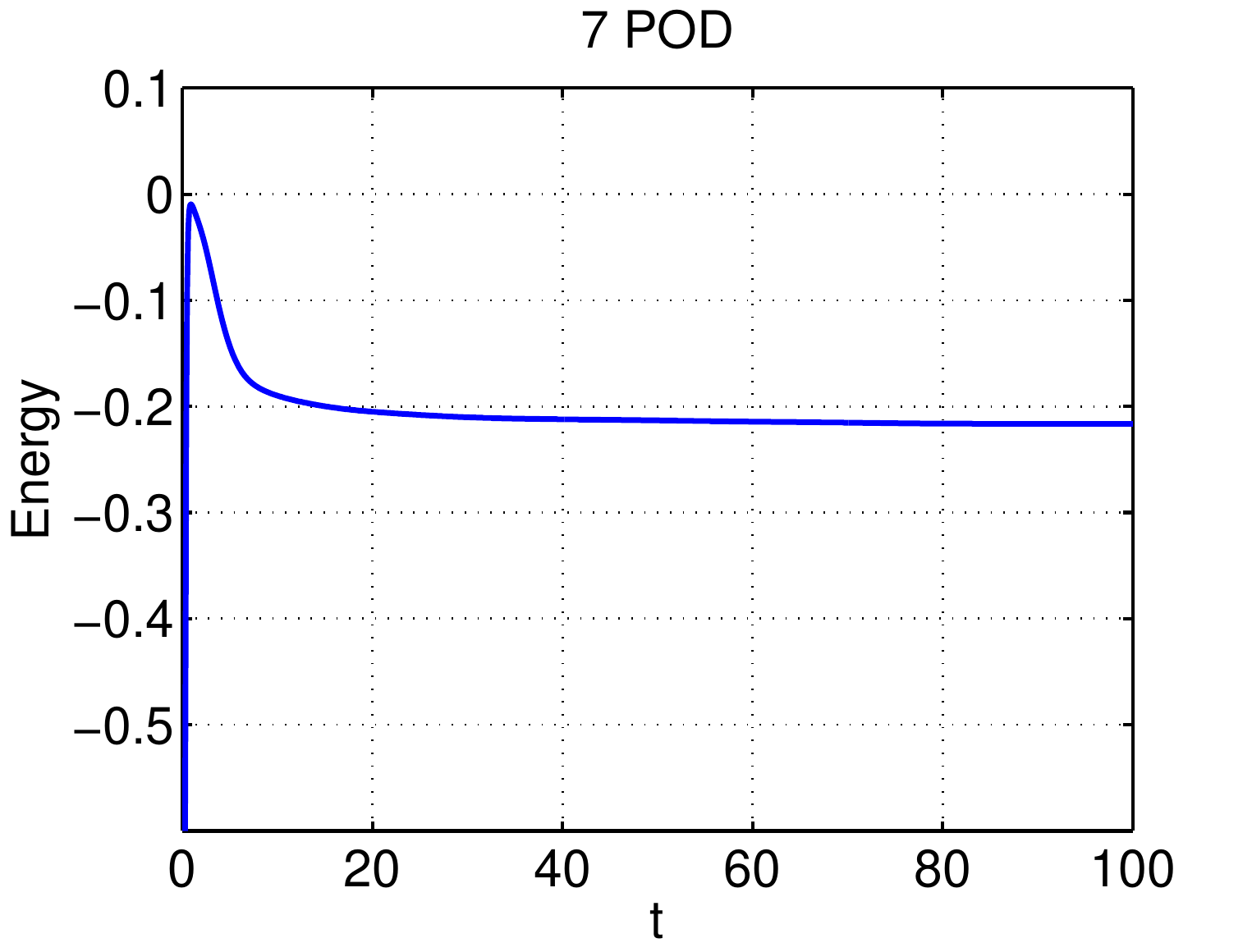}}
\subfloat[]{\includegraphics[scale=0.25]{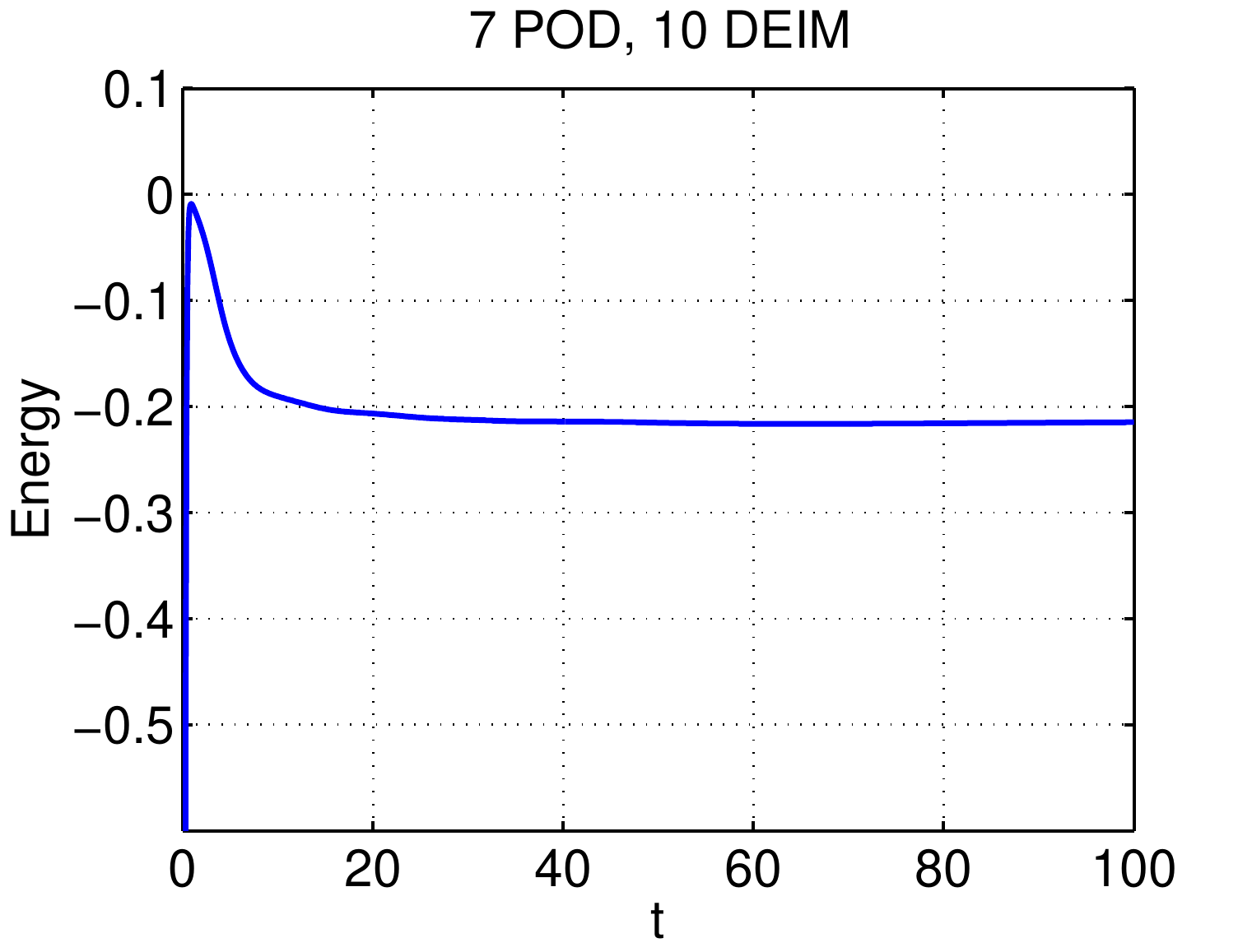}}
\caption{Energy plots\label{sec7_energy}}
\end{figure}

\begin{figure}[htb!]
\centering
\includegraphics[scale=0.35]{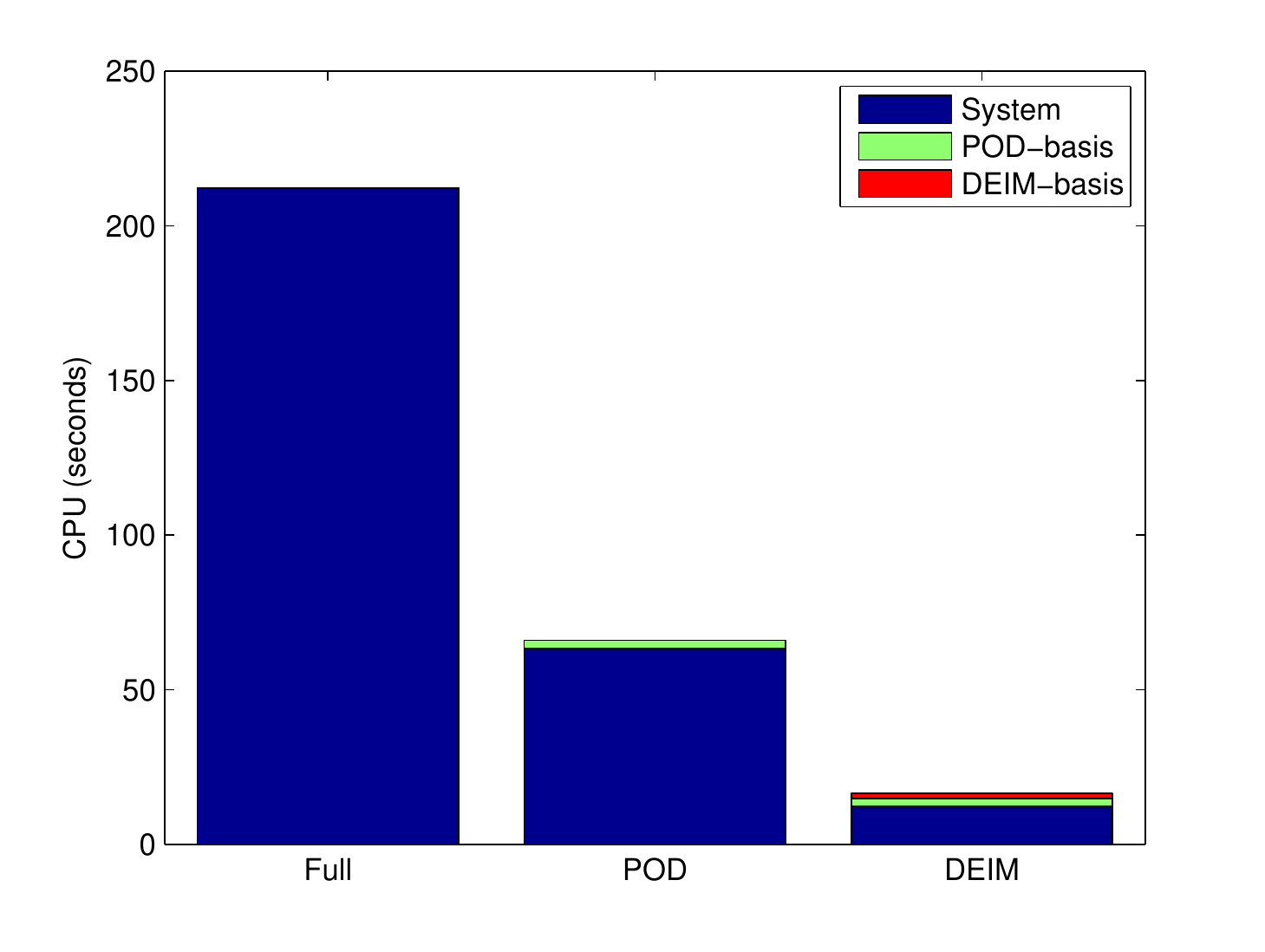}
\caption{CPU times\label{sec7_cpu}}
\end{figure}

In Fig.~\ref{sec7_sing} the decay of the singular values related to $U$, $V$ and the nonlinear term $F$ are shown. In Figs.~\ref{sec7_plotu}, \ref{sec7_plotv}, \ref{sec7_energy}, we demonstrate that the ROM solutions and energy plots for $k=7$ POD and $m=10$ DEIM basis functions are almost the same as the FOM solutions in Fig~\ref{fig9}. The computational efficiency of the POD-DEIM is obvious from Fig.~\ref{sec7_cpu}.

We demonstrate the computational efficiency of the POD-DEIM algorithm for linear and quadratic DG elements. In Table~\ref{table1} we give the CPU time and the speed-up factors $S_{POD}$ and $S_{DEIM}$ of POD and POD-DEIM, respectively,
$$
S_{POD}=\frac{\text{CPU time for FOM}}{\text{CPU time for POD ROM }}, \quad
S_{DEIM}=\frac{\text{CPU time for FOM }}{\text{CPU time for POD-DEIM ROM}}.
$$

In all computation only one Newton iteration is needed.
The results show clearly that as the mesh size increases, the speed-up factors increases, which shows the efficiency of POD-DEIM method. In Table~\ref{table2}, we show the computed mean relative $L^2$ errors between FOM and POD, and between FOM and POD-DEIM solutions which are acceptable because  POD-DEIM  is an approximation of the non-linearity in contrast to the POD.

\begin{table}
\centering
\caption{The computation time (in sec) and speed-up factors for $k=10$ POD basis functions and $m=50$ DEIM basis functions on different grids and with different polynomial degrees}
\label{table1}
\begin{tabular}{ l c c c c c c }\hline
Polynomial degree &  &  $p=1$ & &  & $p=2$ &  \\
\hline
Mesh number & 1 & 2 & 3 & 1 & 2 & 3 \\
\hline
\# triangles & 128 & 512  & 2048   & 128 & 512  & 2048   \\
\# nodes  & 81  &  289  & 1089   & 81  &  289  & 1089   \\
\# DoFs  & 384 &  1536  & 6144   & 768 &  3072  & 12288  \\
\hline
Full        &  9.92 & 26.25 & 124.37 & 29.42  & 81.53 & 423.97   \\
POD         &  7.65 & 13.50  & 33.64  & 22.23  & 38.67  & 116.82  \\
POD-DEIM    &  7.11  & 9.64  & 12.67  & 16.29   & 19.12  &  31.23  \\
$S_{POD}$   &  1.34       & 2.00        & 3.83         &  1.34        & 2.15        &  3.71  \\
$S_{DEIM}$  &  1.40       & 2.90        & 9.81         &  1.86        & 5.50        &  15.39 \\
\hline
\end{tabular}
\end{table} 

\begin{table}
\centering
\caption{Mean relative $L^2$ errors between the solutions of the full system and the ones for POD system and POD-DEIM system with $k=10$ POD basis functions and $m=50$ DEIM basis functions}
\label{table2}
\begin{tabular}{ l c c c c c c }\hline
Polynomial degree &  &  $p=1$ & &  & $p=2$ &  \\
\hline
Mesh number & 1 & 2 & 3 & 1 & 2 & 3 \\
\hline
POD (Comp. u)      &  3.8e-2 & 1.2e-2  & 1.0e-2  &  8.9e-3   &  1.1e-2  &  1.1e-2   \\
POD-DEIM (Comp. u) &  3.6e-2 & 1.5e-2  & 1.5e-2  &  1.2e-2   &  2.5e-2  &  2.9e-2   \\
POD (Comp. v)      &  1.4e-2 & 5.6e-3  & 4.7e-3  &  6.0e-3   &  6.1e-3  &  6.4e-3  \\
POD-DEIM (Comp. v) &  1.4e-2 & 6.0e-3  & 5.4e-3  &  6.6e-3   &  8.1e-3  &  9.1e-3  \\
\hline
\end{tabular}
\end{table} 


\section{Conclusions}
\label{sec:conc}

 Neurons, the building blocks of the central nervous system, are highly complex dynamical systems. The FHN equation as the simplified version of the Hodgkin–-Huxley equation models in a detailed manner activation and deactivation dynamics of a spiking neurons. The FHN equation can describe the bifurcations with the variation of the key parameters for neuron dynamics. Within time FHN equations  became a favorite model for simulation of wave  propagation  in excitable media, such as heart tissue or nerve fiber. Understanding the complex behavior of patterns in neuroscience can have big impact for dealing and preventing with various diseases.
In this paper we have shown that patterns like pulses, fronts, spots and labyrinths for the FHN  equation  are computed accurately using the structure preserving time integrator AVF in combination with DG finite elements in space. MOR can dramatically reduce simulation costs by preserving the behavior of parametrized PDEs, which was demonstrated here for the FHN equation with Turing patterns. In the literature MOR is applied usually in connection with the finite differences and continuous finite elements. We have shown that the DG discretization in space can save the computational cost and due its local structure DG is more efficient than the continuous finite elements.

Quantitative methods are indispensable to structure the clinical, epidemiological, and economic evidence in health care and qualitative insight in making better decisions. Model based analysis allows to quantify  different concepts relating to uncertainty in decision modeling. Therefore effective simulation methods like the MOR can contribute in for more precise decisions and better treatment of diseases. Beside the MOR, the two important fields of the modern scientific computing are uncertainity quantification and optimal control. The FHN equation is an ideal model in this respect to recover the uncertainities with respect to parameters and input variables in the context of neuron modeling. Also optimal treatment of the diseases can be designed using reduced order simulation effectively. Optimization and uncertainty quantification techniques combined with MOR can improve model predictions, to evaluate monitoring schemes and apply better therapy in the  medical practice.

\begin{acknowledgements}
The authors would like to thank the reviewer for the comments and suggestions that help improve the manuscript. This work has been supported by METU BAP-07-05-2015 009.
\end{acknowledgements}



\begin{thebibliography}{18}
\providecommand{\natexlab}[1]{#1}
\providecommand{\url}[1]{{#1}}
\providecommand{\urlprefix}{URL }
\expandafter\ifx\csname urlstyle\endcsname\relax
  \providecommand{\doi}[1]{DOI~\discretionary{}{}{}#1}\else
  \providecommand{\doi}{DOI~\discretionary{}{}{}\begingroup
  \urlstyle{rm}\Url}\fi
\providecommand{\eprint}[2][]{\url{#2}}

\bibitem[{Antil et~al(2014)Antil, Heinkenschloss, and
  Sorensen}]{Heinkenschloss14}
Antil H, Heinkenschloss M, Sorensen C Danny (2014) Application of the discrete
  empirical interpolation method to reduced order modeling of nonlinear and
  parametric systems. In: Quarteroni A, Rozza G (eds) Reduced Order Methods for
  Modeling and Computational Reduction, MS \& A - Modeling, Simulation and
  Applications, vol~9, Springer International Publishing, pp 101--136

\bibitem[{Arnold(1982)}]{arnold82ipf}
Arnold DN (1982) An interior penalty finite element method with discontinuous
  elements. SIAM J Numer Anal 19:724--760

\bibitem[{Barrault et~al(2004)Barrault, Maday, Nguyen, and Patera}]{Barrault04}
Barrault M, Maday Y, Nguyen NC, Patera AT (2004) An ‘empirical
  interpolation’ method: application to efficient reduced-basis
  discretization of partial differential equations. Comptes Rendus Mathematique
  339(9):667--672, \doi{10.1016/j.crma.2004.08.006}

\bibitem[{{Celledoni, E. and Grimm, V. and McLachlan, R. I. and McLaren, D. I.
  and O'Neale, D. J. and Owren, B. and Quispel, G. R.
  W.}(2012)}]{Cellodoni12ped}
{Celledoni, E and Grimm, V and McLachlan, R I and McLaren, D I and O'Neale, D J
  and Owren, B and Quispel, G R W} (2012) {Preserving energy resp. dissipation
  in numerical PDEs using the {"Average Vector Field"} method.} J Comput
  Physics 231:6770--6789

\bibitem[{Chaturantabut and Sorensen(2010)}]{chaturantabut10nmr}
Chaturantabut S, Sorensen DC (2010) Nonlinear model reduction via discrete
  empirical interpolation. SIAM J Scientific Computation 32(5):2737--2764

\bibitem[{Chen and Hu(2014)}]{Chen14sas}
Chen CN, Hu X (2014) {Stability analysis for standing pulse solutions
  to{F}itz{H}ugh-{N}agumo equations}. Calculus of Variations and Partial
  Differential Equations 49:827--845, \doi{10.1007/s00526-013-0601-0}

\bibitem[{Grepl(2012)}]{Grepl3}
Grepl MA (2012) Model order reduction of parametrized nonlinear
  reaction–-diffusion systems. Computers \& Chemical Engineering 43:33 -- 44,
  \doi{j.compchemeng.2012.03.013}

\bibitem[{Hairer and Lubich(2014)}]{Hairer14}
Hairer E, Lubich C (2014) Energy-diminishing integration of gradient systems.
  IMA J Numer Anal 34(2):452--461, \doi{10.1093/imanum/drt031}

\bibitem[{van Heijster and Sandstede(2011)}]{heijster11}
van Heijster P, Sandstede B (2011) Planar radial spots in a three-component
  {FitzHugh–Nagumo} system. Journal of Nonlinear Science 21(5):705--745,
  \doi{10.1007/s00332-011-9098-x}

\bibitem[{van Heijster et~al(2008)van Heijster, Doelman, and
  Kaper}]{heijster08}
van Heijster P, Doelman A, Kaper TJ (2008) Pulse dynamics in a three-component
  system: Stability and bifurcations. Physica D: Nonlinear Phenomena
  237(24):3335 -- 3368, \doi{http://dx.doi.org/10.1016/j.physd.2008.07.014}

\bibitem[{Kunisch and Volkwein(2001)}]{Kunisch01}
Kunisch K, Volkwein S (2001) Galerkin proper orthogonal decomposition methods
  for parabolic problems. Numerische Mathematik 90(1):117--148,
  \doi{10.1007/s002110100282}

\bibitem[{Marquez-Lago and Padilla(2014)}]{lago14}
Marquez-Lago TT, Padilla P (2014) A selection criterion for patterns in
  reaction–diffusion systems. Theoretical Biology and Medical Modelling 11:7,
  \doi{10.1186/1742-4682-11-7}

\bibitem[{Or-Guil et~al(1998)Or-Guil, Bode, Schenk, and Purwins}]{guil98}
Or-Guil M, Bode M, Schenk CP, Purwins HG (1998) Spot bifurcations in
  three-component reaction-diffusion systems: The onset of propagation. Phys
  Rev E 57:6432--6437, \doi{10.1103/PhysRevE.57.6432}

\bibitem[{Riviere(2008)}]{riviere08dgm}
Riviere B (2008) {Discontinuous Galerkin Methods for Solving Elliptic and
  Parabolic Equations: Theory and Implementation}. SIAM

\bibitem[{Tiso and Rixen(2013)}]{Tiso13}
Tiso P, Rixen DJ (2013) Discrete empirical interpolation method for finite
  element structural dynamics. In: Topics in Nonlinear Dynamics, Volume 1
  Proceedings of the 31st IMAC, A Conference on Structural Dynamics, Topics in
  nonlinear dynamics, Vol. 1, The Society for Experimental Mechanics, pp
  203--212

\bibitem[{Turing(1952)}]{Turing52}
Turing AM (1952) The chemical basis of morphogenesis. Philosophical
  Transactions of the Royal Society of London Series B, Biological Sciences
  237(641):37--72, \doi{10.1098/rstb.1952.0012}

\bibitem[{Yanagida(2002{\natexlab{a}})}]{Yanigada02mfrd}
Yanagida E (2002{\natexlab{a}}) {Mini-Maximizers for Reaction-Diffusion Systems
  with Skew-Gradient Structure}. Journal of Differential Equations
  179:311--335, \doi{10.1006/jdeq.2001.4028}

\bibitem[{Yanagida(2002{\natexlab{b}})}]{Yanigada02sps}
Yanagida E (2002{\natexlab{b}}) {Standing Pulse Solutions in Reaction-Diffusion
  Systems with Skew-Gradient Structure}. Journal of Dynamics and Differential
  Equations 14:189--205, \doi{10.1023/A:1012915411490}

\end{thebibliography}
\end{document}